\documentclass[final,1p,times]{elsarticle}
\usepackage{lineno,hyperref}
\usepackage{amsmath,amsfonts}
\usepackage{verbatim}
\usepackage{latexsym}
\usepackage{graphicx}
\usepackage{epstopdf}
\usepackage{subfig}
\usepackage{caption}
\usepackage{dsfont}
\usepackage{indentfirst}
\usepackage{xcolor}
\usepackage{algorithm}
\usepackage{algpseudocode}
\usepackage{booktabs}
\usepackage{natbib}
\usepackage{color}
\usepackage{multirow}  
\usepackage{textcomp}%
\usepackage{manyfoot}%
\usepackage{booktabs}%
\usepackage{threeparttable}
\usepackage{tablefootnote} 
\allowdisplaybreaks[4]  
\newtheorem{theorem}{Theorem}[section]
\newtheorem{lemma}[theorem]{Lemma}
\newtheorem{corollary}[theorem]{Corollary}

\newtheorem{assumption}[theorem]{Assumption}
\newtheorem{example}[theorem]{Example}
\newtheorem{proposition}[theorem]{Proposition}
\newtheorem{remark}[theorem]{Remark}

\def\lan{\langle} \def\ran{\rangle}

\newcommand{\D}{\Delta}

\def\eps{\varepsilon}

\def\th{\theta}
\def\ve{\vee}
\def\we{\wedge}
\def\a{\alpha} \def\g{\gamma}
\def\b{\beta}

\def\k{\kappa} \def\l{\lambda}

\def\K{\times}

\def\f{\forall}

\newcommand{\Ito}{It\^{o} }
\newcommand{\Holder}{H\"{o}lder }
\newcommand{\BDG}{Burkholder-Davis-Gundy }

\newcommand{\E}{\mathbb{E}}

\newcommand{\PP}{\mathbb{P}}
\newcommand{\II}{\mathbf{1}} 

\newcommand{\R}{\mathbb{R}}

\newcommand{\fd}{f_{\Delta}}
\newcommand{\gd}{g_{\Delta}}
\newcommand{\ggd}{g'\cdot g_{\Delta}} 
\newcommand{\tkk}{t_{k+1}}
\newcommand{\tk}{t_{k}}
\newcommand{\qu}{\quad}
\newcommand{\no}{\nonumber }

\newcommand{\pd}{\pi_{\Delta}}

\makeatletter
\@addtoreset{equation}{section}
\makeatother
\renewcommand{\thesection}{\arabic{section}}
\renewcommand{\theequation}{\thesection.\arabic{equation}}

\modulolinenumbers[5]

\begin{document}
\begin{frontmatter}
\title{Positivity-preserving truncated Euler and Milstein methods for financial SDEs with super-linear  coefficients}


\author[mainaddress]{Shounian Deng}
\author[address]{Chen Fei }

\author[mainaddress]{Weiyin Fei \corref{correspondingauthor}}
\cortext[correspondingauthor]{Corresponding author}
\ead{wyfei@ahpu.edu.cn}
\author[thirdaryaddress]{Xuerong Mao}
\address[mainaddress]
{School of Mathematics-Physics and Finance, Anhui Polytechnic University,  Wuhu  241000, China}
\address[address]{Business School, University of Shanghai for Science and Technology, Shanghai 200093, China}
\address[thirdaryaddress]{Department of Mathematics and Statistics, University of Strathclyde, Glasgow G1 1XH, UK}
\begin{abstract}

In this paper, we propose two variants of the positivity-preserving  schemes, namely the truncated Euler-Maruyama (EM)  method and the truncated Milstein scheme, applied to stochastic differential equations (SDEs) with positive solutions and super-linear coefficients.  Under some regularity and integrability assumptions we derive  the optimal strong convergence rates of the two schemes. Moreover, we demonstrate flexibility of our approaches by applying the truncated methods to   approximate SDEs with super-linear coefficients (3/2 and  A\"it-Sahalia models) directly  and also   with sub-linear coefficients (CIR model) indirectly.
Numerical experiments
are provided to verify the effectiveness of the theoretical results.

\end{abstract}

\begin{keyword}
Truncated  EM method, Truncated Milstein method,  Strong convergence order, Positivity preservation.
\end{keyword}

\end{frontmatter}

\section{Introduction}

The goal of this paper is to derive a positivity-preserving numerical method for
scalar SDEs which takes values in $\R_+$ and have super-linear coefficients.
Typical examples of such SDEs in mathematical finance and bio-mathematical applications
are the Heston-3/2 volatility process  \cite{Heston1997ASN}
\begin{align}\label{model1}
dX(t) = c_1 X(t) (c_2 - X(t)  ) dt +  \sigma X(t)^{3/2} dB(t), \qu X(0) = X_0 >0,
\end{align}
with $c_1, c_2, \sigma >0$,
and the A\"it-Sahalia (AIT, for short) model \cite{As1996}
\begin{align}\label{Ait1}
dX(t) = & \left  (a_{-1}{X(t)^{-1}} - a_0 + a_1 X(t) - a_2 X(t)^{\kappa} \right  ) dt + b  X(t) ^{\theta } dB(t), \qu X(0) =  X_0>0,
\end{align}
where all constant parameters are nonnegative,  $\kappa, \theta >1$, and $B(t)$ is a Wiener process.
The SDEs appearing in such models are highly nonlinear and may contain the singularity in the neighbourhood of zero   in the drift or diffusion coefficient.
Such SDEs in almost all cases cannot be solved explicitly, and it has been and still is a very 
active topic of research to approximate SDEs with super-linear coefficients; 
see, e.g., \cite{Mao2011BIT,Lucas2014NM,FCJ2016SIAM,HJ2019FS,Gan_Wang2020NA, CS2020IMA,Mao_Wei2021CAM, Hong2021BIT, Hu_Zhao2021CNSNS,CS2022APNUM,Gan2023CAM,Lord2023ap,Tian_jcp2023,Xiao2024cam,Cai_Mao2024cns,WXJ2024na_Mil,WXJ2024na2,Cai2024CAM},
and the references therein.
 \par
\begin{table}[!t]
  \centering
   \caption{ Summary of the condition over the parameters for the convergence of the positivity-preserving schemes for AIT }
   \label{Tab2}
\begin{tabular}{c|c|c|c}
  \hline \hline
  Scheme & Norm  &  Parameter Range & Rate \\ \hline \hline
  Lamperti  Projected EM  & \multirow{2}{*}{ $L^p, p \ge 1$} & \multirow{6}{*}{$ \k + 1  > 2 \th $}
       & \multirow{2}{*}{$\frac{1}{p}$}    \\
(Chassagneux et al. \cite{FCJ2016SIAM})    &  &  \\  \cline{1-2} \cline{4-4}
  Lamperti  Semi-Discrete EM  & \multirow{2}{*}{$L^2$}  &
       & \multirow{2}{*}{$\frac{1}{2}$}    \\
 (Halidias and Stamatiou \cite{HN2023dcds})    &  &  \\  \cline{1-2}\cline{4-4}
Lamperti BEM   & \multirow{8}{*}{$L^p,p \ge 2$}    &
       &   \multirow{2}{*}{$ 1 $ } \\
 (Neuenkirch and Szpruch \cite{Lucas2014NM})   &  &   \\  \cline{1-1} \cline{3-4}
 Logarithmic  Truncated EM  &   & \multirow{2}{*}{$ \k + 3  > 4 \th $} & \multirow{4}{*}{$\frac{1}{2}$}     \\
(Yi et al. \cite{Hu_Zhao2021CNSNS}, Lei et al. \cite{Gan2023CAM})    &  &  \\  \cline{1-1} \cline{3-3}
     Logarithmic Truncated EM  &    &     &    \\
    (Tang and Mao \cite{Tang_Mao2024ap})    &     & \multirow{16}{*}{$ \k + 1  > 2 \th $}    &     \\  \cline{1-1}\cline{4-4}
  BEM  &  &
       & \multirow{2}{*}{without rate}    \\
(Szpruch et al. \cite{Mao2011BIT})  &  &  \\  \cline{1-2} \cline{4-4}
  Theta EM  & \multirow{2}{*}{ $L^2 $} &
       & \multirow{2}{*}{$\frac{1}{2}$}    \\
(Wang et al. \cite{WangXJ2020BIT})  &  &  \\  \cline{1-2} \cline{4-4}
  Truncated EM  & \multirow{4}{*}{ $L^p, p \ge 1 $} &
       & \multirow{2}{*}{$\frac{1}{2p}$}    \\
(Deng et al. \cite{Deng2023BIT})  &  &  \\  \cline{1-1} \cline{4-4}
\multirow{2}{*}{ Proposed Truncated EM } &   &    &   \multirow{4}{*} {$\frac{1}{2}$} \\
                                & & &                                \\ \cline{1-2}
Semi-implicit  Tamed EM  & \multirow{6}{*}{ $L^2 $} &       &     \\
(Liu et al. \cite{WXJ2024na2})    &  &  \\  \cline{1-1} \cline{4-4}
Semi-implicit  Projected Milstein  &  &
       & \multirow{4}{*}{$1$}    \\
(Jiang et al. \cite{WXJ2024na_Mil})    &  &  \\  \cline{1-1}
\multirow{2}{*}{Proposed  Truncated Milstein}  &  &
       &    \\
&   &  &  \\  \cline{1-4}

\end{tabular}
\end{table}
Two significant challenges in constructing a numerical method for non-linear
 SDEs with positive solutions  are to preserve positivity and to derive convergence with as high a rate as possible.
 In this section,
 we mainly discuss the  positivity-preserving numerical methods for AIT \eqref{Ait1}.
 In fact, the techniques to handle the superlinearity and  the singularity of the AIT coefficients   can be applied to a more general SDE with positive solutions and super-linear coefficients.
 Table \ref{Tab2} gives a
summary of  the known results
 of some of the key methods which  possess positivity preservation   for AIT \eqref{Ait1}. We now discuss two main methods for the strong approximation of \eqref{Ait1}, the first is based on  Lamperti/logarithmic transformation, and the second is direct approximation of AIT \eqref{Ait1}.
 \par
 A classical transformation approach  is to apply the Lamperti transformation $Y = X^{1-\th}$
 in order to obtain an auxiliary SDE in $Y$ with a globally Lipschitz diffusion, but
 a drift function which is unbounded when solutions are in a neighbourhood of zero.
 For  Lamperti transformation methods that preserve positivity for \eqref{Ait1} see, for example, Lamperti backward EM (BEM) \cite{Lucas2014NM},   Lamperti  projected EM \cite{FCJ2016SIAM},  and Lamperti Semi-Discrete EM  \cite{HN2023dcds}.
 However, this type of polynomial transformation would shift all the nonsmoothness into the drift function.
 In recent years, much effort has focused on deriving the combination of logarithmic transformation and Euler/Milstein-type schemes for general SDEs which have positive solutions.
 For such SDEs, an effective explicit scheme   that preserves positivity is the logarithmic truncated EM proposed by Yi et al. \cite{Hu_Zhao2021CNSNS}.  The method is to  apply a logarithmic transformation to original SDE and to numerically approximate the transformed SDE by  a truncated EM (an explicit method established in  \cite {Mao2015truncated,Mao2016rates}).
  Under the  Khasminskii-type and the global monotonicity conditions on coefficients of  transformed SDE, the method was shown to have strong $L^p$-order of convergence $1/2$.
 In a subsequent paper \cite{Xiao2024cam}, a logarithmic truncated Milstein scheme was proved to have strong $L^p$-convergence of order $1$,  see Table \ref{Tab2}.
 However, when the  logarithmic  method in  \cite{Hu_Zhao2021CNSNS} is applied to  AIT \eqref{Ait1},
the method   requires a more restrictive condition:
 $ \k  > 4 \th -3$  in order to satisfy the Khasminskii-type requirement for the transformed SDE, in contrast to  the other positivity-preserving methods in Table \ref{Tab2}.
 The same approach was also used by Tang and Mao in \cite{Tang_Mao2024ap} to derive strong
 order  $1/2$ for the SDEs with positive solutions, but weaker assumptions on the coefficients of  original SDE was required.
%
%
 \par
 The second main approach that preserve positivity of the numerical solution is to directly discretise  \eqref{Ait1} by using the implicit or truncated tricks.
 A backward EM (BEM)    was
 shown to preserve  positivity of the solution for AIT \eqref{Ait1}
 in   Szpruch et al. \cite{Mao2011NA_Sahalia}, where the authors
 proved the strong convergence of the  method, but without revealing a rate.
 This gap was filled by Wang et al. in \cite{WangXJ2020BIT},  where the authors
successfully recovered the mean-square convergence rate of order $1/2$ for the theta EM applied to \eqref{Ait1}.  However, the computational costs of
solving implicit algebraic equations produced by classical BEM rise as the parameter $\k$
increases.
\par
A truncated EM that seeks to directly approximate AIT was proposed by Emmanuel and Mao in \cite{Coffie_Mao2020CAM},  where the authors proved  the strong $L^p $-convergence of the method, but they did not reveal a   convergence rate.
Recently, the $L^p $-convergence rate of order $\frac{1}{2p}$ of this method was obtained by  Deng et al.  \cite{Deng2023BIT}.
At the same time,  a novel semi-implicit plus tamed EM was introduced by Liu et al. in \cite{WXJ2024na2},
 where the authors established the strong $L^p$-convergence  of rate $1/2$ based on an implicit treatment of the singular term $a_{-1} x^{-1}$ to ensure the property of positiveness and   a suitable taming    of the  super-linear terms $-a_2 x^{\k}$ and $b x ^{\th}$.
 It should be pointed out that this method is  explicit, as the proposed scheme  can be explicitly expressed  by finding a positive root of a quadratic equation. Therefore, compared with the implicit schemes in \cite{Mao2011BIT,Lucas2014NM,WangXJ2020BIT} the computational costs in \cite{WXJ2024na2}  can be significantly reduced. Similar approach was also used by Jiang et al. in \cite{WXJ2024na_Mil} to derive strong order $1$ of semi-implicit projected Milstein for AIT process. See Table \ref{Tab2} for a comparison.
\par
The motivation of this paper is the following observation:
in  \cite{Deng2023BIT}
 the theoretical $L^p$-convergence rate of truncated EM for AIT is only $ \frac{1}{2p}$,  rather than the optimal $ \frac{1}{2}$, and
 it  should be further increased. As a consequence,  this paper is concerned with
the optimal strong  convergence rate of the positivity-preserving truncated methods for general SDE \eqref{Eq0}  with solution and super-linear coefficients.
We note that the error analysis in \cite{Deng2023BIT} is based on the higher moment estimates of the numerical solution and the interpolation of the discretization scheme to continuous time. According to the theory of Mao \cite{Mao2015truncated,Mao2016rates}, it is inevitable to estimate the  probability of truncated EM solution escaping from the truncated interval $(1/R, R)$, i.e., $\PP ( \rho_{\D, R} \le T) $.
However,  \cite[Lemma 3.2]{Deng2023BIT} implies that  the upper bound for this probability is  of order $1/2$ at most,  i.e., \begin{align}\label{Temp731}
\PP ( \rho_{\D, R} \le T) \le CR^{(2 \g +1) \ve (\g +4)} \D^{1/2},
\end{align}
which prevents the truncated EM in \cite{Deng2023BIT} from  achieving the optimal strong order $1/2$.
In this work, we seek to bypass this issue.
\par
Motivated by the approaches of \cite{FCJ2016SIAM,BIK2016PEM}, we adopt the  error analysis based on the
higher moment  and inverse moment estimates of the true solution  and the discretised    scheme rather than the continuous time extension of numerical solution.  We first prove a key Lemma \ref{Lem1} which shows
the global monotonicity  for the truncated functions and implies the stochastic C-Stability of the truncated EM in the sense of \cite[Definition 3.2]{BIK2016PEM}.
By the finiteness of the $p_0$th moments and  $p_1$th inverse moments of the true solution $X$,
we then obtain the local truncations errors of $f(X)$ and $g(X)$   in $L^p$, which implies the stochastic B-Consistency of the truncated EM in the sense of \cite[Definition 3.3]{BIK2016PEM}, see Lemma \ref{Lem3.5}.
Combining Lemma   \ref{Lem1} with Lemma  \ref{Lem3.5} allows us to establish the optimal $L^p$-convergence of order $1/2$ of positivity-preserving truncated EM for SDE \eqref{Eq0} in the similar proofs as in
\cite[Theorem 5.36]{Shi2024jde}.
According to the global error estimate \eqref{Result2} in Theorem \ref{Theorem3.1}, we show that
the truncated numerical solution $Y_\D$ has some finite moments and  inverse moments, see Proposition \ref{Cor3.1}.
By $L^p$-convergence Theorem \ref{Theorem3.1}  and  moment boundedness Proposition \ref{Cor3.1}, we apply truncated EM
to directly approximate the 3/2 and the AIT models, see Applications \ref{32model} and \ref{AIT_app}.
Based on the Lamperti transformation, we also apply this method to  indirectly approximate the CIR process, see Application  \ref{cir_ap}.
Combining the Milstein scheme  with the truncated strategy, we  also prove the strong $L^2$-convergence  order $1$ for the  truncated Milstein, in the similar method of analyzing the convergence rate of the truncated EM applied to SDE \eqref{Eq0}.
\par
The main contributions of this paper are as follows:
\begin{itemize}
  \item In contrast with the results of \cite{Deng2023BIT}, our truncated EM for AIT obtains  the optimal $L^p$-convergence of order $\frac{1}{2}$, which is  far better  than the rate $\frac{1}{2p}$   in \cite{Deng2023BIT}.
  \item Compared with the  projected EM in \cite{FCJ2016SIAM}, where the method requires that the diffusion coefficients of SDEs are globally Lipschitz continuous, our truncated methods can apply to a wide class of  SDEs with positive solutions and super-linearly growing diffusion coefficients.
\item
  Our error analysis method does not relies on  the  high-order moments and  inverse moments estimates, the continuous time extension  of the numerical solution. Thus our truncated methods can be applied to  approximate  SDEs with  super-linear coefficients directly  and also with  sub-linear coefficients indirectly, see Applications \ref{Sec4}.
\end{itemize}

\par
The structure of this paper is as follows. In Section \ref{Sec2} we give the form of the SDE with positive solution and super-linear coefficients, specify   some regularity and integrability constraints placed upon the coefficients for the main strong convergence theorems. In Section \ref{Sec3} we set up the framework  and  prove the optimal  strong convergence order as well as  establish the  moment boundedness results   for  the positivity-preserving truncated EM and Milstein schemes. In Section \ref{Sec4} we apply the proposed methods
to some financial SDEs, such as the 3/2 and the AIT  as well as  the CIR models.  In Section \ref{Sec5} we numerically compare convergence and efficiency of several commonly used methods.  Finally, we have included in a small appendix a couple of proofs to make this paper self contained.
%

%
\section{Setting and preliminaries}

\label{Sec2}

Let $(\Omega , {\mathcal F}, \mathbb{P})$ be a complete probability
space with a filtration  $\{{\cal F}_t\}_{t\ge 0}$ satisfying the usual conditions (i.e., it is increasing and right continuous while $\cal{F}_\textrm{0} $ contains all $\mathbb{P}$\textrm{-}null sets).
Let $\R = ( -\infty, \infty)$ and $\R_+ = (0, \infty)$.
We write $| \cdot |$ and $\lan \cdot, \cdot \ran$ respectively for the Euclidean norm and the inner produce on $\R$.
Let $\{B(t)\}_{t \ge 0}$  be a one-dimensional Brownian motion with respect to the normal
filtration $\{{\cal F}_t\}_{t\ge 0}$.
Let $\E$ denote the expectation. Denote by $L^p(\Omega; \R)$, for $p >0$, the set of random variables $Z$ such that $\| Z \|_p: = ( \E [|Z|^p])^{1/p} < \infty$.
We denote by $\E_{t} [X] := \E [X | \mathcal F_{t}]$ the conditional expectation given the filtration $ \mathcal F_{t}$.
For two real numbers $a$ and $b$,  $a \ve b:  = \max (a,b)$ and $a \we b : = \min (a,b)$.
W denote by $\mathcal{C}^{2}(\R_+)$  the space of twice differential functions with continuous derivatives on $\R_+$.
For any $x \in \R_+$, we denote $\lceil x \rceil$ the rounded up integer.
For a set $G$, its indicator function is denoted by $\II_G$. In the following it will be convenient to introduce the abbreviation
\begin{align*}
f'(x) : = \frac{df(x)}{dx}, \qu  f''(x) := \frac{d^2f(x)}{dx^2},
\qu \textrm{and} \qu  g'\cdot g (x) : = \frac{d g(x)}{dx } g(x), \qu \f x \in \R.
\end{align*}
Moreover, given $V(x) \in \mathcal{C}^{2}(\R, \R_+) $,  we define a functional $ \mathbb{L} V: \R \to \R $ by
\begin{align*}
\mathbb L V(x) = V'(x) f(x) + \frac{1}{2}V''(x) |g(x)|^2,\qu \f x \in \R.
\end{align*}
In this work, we frequently apply the weighted Young inequality
\begin{align}\label{Inequality1}
x^a y^b  \le \eps x^{a+b} + \frac{a}{a+b} \left (\frac{a}{ \eps (a+b)} \right )^{a/b} y^{a+b}, \qu \f x, y, a,b, \eps >0.
\end{align}


Consider  a one-dimensional SDE of the form
\begin{align}\label{Eq0}
dX(t) & = f(X(t))dt + g(X(t))dB(t),  \; t > 0 \\
X(0) & = X_0>0. \no
\end{align}
Throughout this paper, we assume that SDE \eqref{Eq0} has a unique strong solution in $\R_+$.
We now formulate the conditions on the drift and the diffusion coefficient functions.
\begin{assumption}\label{H0}
There are $\a \ge 0$,
$\b \ge  0 $, $K_1 >0$, such that
\begin{align}\label{Eq32}
|f(x)-f(y)| \ve |g(x) - g(y)| \le K_1 \Big ( 1+ x^\a+ y^\a + x^{-\b} + y^{-\b}  \Big)|x-y|, \qu \f x, y \in \R_+.
\end{align}
\end{assumption}
\begin{assumption} \label{Assu1}  Assume that $f, g \in  \mathcal{C }(\R_+)$ and
there are  positive numbers $K$ and $q_0 >1 $, such that
\begin{align}\label{Eq3}
2f'(x) + q_0 |g'(x)|^2 \le K, \quad \f x \in \R_+.
\end{align}
\end{assumption}
 \begin{assumption}\label{H1}
 There are constants $p_0 >  2 (\a +1)$ and $p_1 >2 \b $, such that
 \begin{align}\label{Eq61}
\sup_{0 \le t \le T} \E \left [ |X(t)|^{p_0} + |X(t)|^{-p_1} \right ] < \infty.
 \end{align}
 \end{assumption}
%
%

In what follows, $C$ will be used to denote a  positive constant depending on $K$, $T$, $\a$, $\b$, $p_0$, $p_1$ and $X_0$, but whose value  may be different from line to line. We denote it by $C_p$ if it depends on an extra parameter $p$.

Assumption \ref{Assu1} implies that the monotonicity condition is fulfilled on $\R_+$, i.e.,
\begin{align}\label{Eq4}
2 \lan x-y, f(x) - f(y)\ran  + q_0 |g(x)-g(y)|^2 \le K|x-y|^2, \quad \f x,y \in \R_+,
\end{align}
and
\begin{align}\label{Eq39}
|f(x)| \ve |g(x)| \le \Big (12K_1 \ve |f(1)|\Big ) \Big (1 + x^{1+\a} + x^{-\b} \Big ), \qu \f x \in \R_+,
\end{align}
see \cite[Remark 2.1]{Tang_Mao2024ap}.
Assumption \ref{H1} imposes a condition on the moments of the true solution $X$ in terms of the locally Lipschitz exponents $\a$, $\b$. Combining this  assumption and Assumption \ref{H0} allows us to bound the local truncations errors for $f(X)$ and $g(X)$, see Lemma \ref{Lem3.5}.
\par
Let $\D \in (0,1]$,  define a truncation mapping $
\pd: \R \to \left [\frac{1}{R(\D)}, R(\D) \right ]$ by
\begin{align}\label{eq23}
\pd(x) =  \frac{1}{R(\D)}  \ve (x \we R(\D) ),    \quad  \f x \in \R,
\end{align}
where $R(\D) := L_1 \D^{-\g}$ is the size of truncated interval with
 $L_1 \in \left [ {\frac{1}{X_0}} \ve X_0 , \infty \right )$  and  $  \g \in \left (0, \frac{1}{2 (\a \ve \b )} \right ]$  is a constant to be determined later.
Define the truncated functions
\begin{align}\label{eq2}
f_\D(x) = f(\pd (x) ) \qu \textrm{and} \qu  g_\D(x) = g(\pd (x) ), \qu \f x \in \R.
\end{align}

The following lemma shows that the truncated mapping $\pd$ is $1$-Lipschitz continuous on $\R$, and the truncated functions $\fd$ and $\gd$ preserve the monotonicity condition \eqref{Eq4}.
\begin{lemma}\label{Lem1}
Let  Assumption \ref{Assu1} hold.
Let $\pd$, $\fd $ and $\gd$ be the truncated mapping and  functions defined in \eqref{eq23}, \eqref{eq2} respectively.
Then
\begin{align}
|\pd(x)- \pd(y)| & \le |x-y|, \quad \f x,y \in \R, \label{Eq2}\\
2 \lan x-y, \fd(x) - \fd(y)\ran  + & q_0 |\gd(x)-\gd(y)|^2  \le 2K|x-y|^2, \quad \f x,y \in \R. \label{Eq11}
\end{align}
\end{lemma}
The proof of this lemma is postponed to Appendix \ref{Append1}.
From Lemma \ref{Lem1}, we have,
\begin{align}
|\pd(x) & \le |x| + \D,\qu \f x \in \R,  \label{Eq42} \\
2\lan x, \fd (x) \ran + (q_0 -1) &  |\gd (x)|^2  \le K_2 (1 + |x|^2), \qu \f x \in \R. \label{Eq42_2}
\end{align}
Moreover, by Assumption \ref{H0} and \eqref{eq23} as well as \eqref{Eq2}, we have
\begin{align}\label{eq62}
|\fd (x) - \fd (y)|\ve |\gd (x) - \gd (y)| \le \varphi (\D) |x -y|, \qu \f x,y \in \R,
\end{align}
where $\varphi(\D) = C \D^{- \g (\a \ve \b)}$  with
\begin{align}\label{Cond1}
\varphi (\D)^2 \D \le C.
\end{align}

We now provide two lemmas, their proofs are postponed to  Appendix \ref{Append2} and \ref{Append3}.
%
%
 Lemma \ref{Lem3.3_2} reveals the \Holder continuity of the true solution to \eqref{Eq0} with respect  to the norm  in $L^p(\Omega; \R)$, while Lemma \ref{Lem3.1} shows the $L^p(\Omega; \R)$-error due to truncating the true solution $X$ on $[1/R,R]$.
\begin{lemma}\label{Lem3.3_2}
Let Assumptions \ref{H0}, \ref{Assu1} and \ref{H1} hold.
Let $p \in \left (0,
\frac{p_0}{1+\a} \we  \frac{p_1}{\b }\right ]$.
 Then,
\begin{align}\label{Eq152}
\E \left [ |X(t+\D) - X(t)|^p  \right ]\le C \D^{p/2} , \qu t \in [0, T],
\end{align}
where $C$ is a positive constant independent of $\D$.
\end{lemma}
%
%
\begin{lemma}\label{Lem3.1}
Under the same assumption as Lemma \ref{Lem3.3_2}, it holds that
\begin{align}
& \E \left [ | f(X(t)) - \fd (X(t)|^p \right ] \ve
\E \left [ | g(X(t)) - \gd (X(t)|^p \right ] \le \frac{1}{\psi(\D)}=C \Big (  \D^{\g(p_0 - p\a -p)} + \D^{\g(p_1 - p\b +p)} \Big ), \qu t \in [0, T] , \label{Eq94} \\
& \E \left [ | X(t) - \pd (X(t))|^p \right ] \le
\left ( \frac{C_{p_0}}{R(\D)^{p_0-p}} +  \frac{C_{p_1}}{R(\D)^{p_1 + p}}  \right )
=C \Big ( \D^{\g (p_0 -p)} + \D^{\g (p_1 + p)} \Big ), \qu t \in [0, T] ,\label{Eq94_3}
\end{align}
where $\psi(\D)= C\Big ( \D^{-\g (p_1 - p \b +p)} + \D^{-\g (p_0 - p \a -p)} \Big ) $ with $p_0$, $p_1$ given by Assumption \ref{H1}.
\end{lemma}

\section{Main Results} \label{Sec3}
\subsection{Positivity-preserving Truncated EM (TEM) scheme}
The first  positivity-preserving scheme is termed truncated EM  scheme. To be more precise,
let  $\D \in (0, 1]$  be a step size,
define the following TEM scheme by setting $X_\D  (t_0)   = X_0$ and computing
\begin{align}\label{PTEM}
  Y_\D (\tk )  & =  \pd ( X_\D (\tk )),  \no  \\
 X_\D (\tkk)  & = X_\D (\tk ) + f (Y_\D (\tk ) ) + g (Y_\D (\tk ) ) \D B_k,
\end{align}
for $k = 0,1,2, \cdots$,  where $\D B_k = B(\tkk) - B(\tk)$, $t_k = k \D$ and $\pd $ has been defined in \eqref{eq23}.
 By \eqref{eq2}, \eqref{PTEM} can be rewritten as the following form:
\begin{align}\label{PTEM2}
X_\D (\tkk)  & =X_\D (\tk ) + \fd (X_\D (\tk ) ) + \gd (X_\D (\tk ) ) \D B_k, \qu k=0,1,2, \cdots.
\end{align}

%
%

We now begin to bound the local truncations errors  for $f(X)$ and $g(X)$ in $L^p(\Omega; \R)$ sense.
\begin{lemma}\label{Lem3.5}
Let Assumptions \ref{H0}, \ref{Assu1} and \ref{H1} hold with
\begin{align}\label{Cond2}
p_0 \ge  (4\a +2) \qu \textrm{and} \qu  p_1 \ge 4 \b .
\end{align}
%
Let
$p \in \left [2,  \frac{p_0}{2\a+1} \we \frac{p_1}{2\b} \right ]$. Then
\begin{align}
\E \left [ |M^{(d)}_k|^p \right ]  & \le C\D^p  \Big (\frac{1}{\psi (\D )} +  \D ^{p/2}\Big ) , \quad  k = 0,1, \cdots,\label{Eq148} \\
\E \left [ |M^{(w)}_k|^p \right ] & \le C\D^{p/2}  \Big (\frac{1}{\psi (\D )} +  \D^{p/2} \Big ) , \quad k = 0,1, \cdots, \label{Eq148_2}
\end{align}
where
\begin{align*}
& M^{(d)}_k: = \int_{\tk}^{\tkk}\Big ( f(X(s)) - \fd (X(\tk)) \Big )ds, \\
& M^{(w)}_k: = \int_{\tk}^{\tkk}\Big ( g(X(s)) - \gd (X(\tk)) \Big )dB(s),
\end{align*}
and $\psi $ is defined  \eqref{Eq94}.
\end{lemma}
\textbf{Proof.}
One notes that
condition \eqref{Cond3} implies
\begin{align}\label{Tempp616}
\frac{p}{p_0 - p \a} \le \frac{1}{\a +1} \quad \textrm{and} \quad
  \frac{p}{p_1 - p \b} \le \frac{1}{\b } .
\end{align}
For any $ \tk < s \le \tkk$, we conclude from \eqref{Eq32}   and the \Holder inequality  that
\begin{align}\label{Eq162_6}
 & \E | f(X(s)) - f(X(\tk))|^p \ve \E | g(X(s)) - g(X(\tk))|^p \\
 & \le C \E \left [ ( 1 + |X(s)|^{p\a} + |X(\tk) |^{p\a} +  |X(s)|^{-p\b} + |X(\tk) |^{-p\b} )|X(s) - X(\tk)|^p \right ] \no \\
 & \le C  \left ( 1 +  \ (\E |X(s)|^{p_0})^{p\a/p_0} +  ( \E |X(\tk)|^{p_0} )^{p\a/p_0} \right ) \left ( \E |X(s) - X(\tk) |^{pp_0/(p_0 - p \a)} \right )^{(p_0 - p\a)/p_0} \no \\
 & \qu + C  \left (   \ (\E |X(s)|^{-p_1})^{p\b/p_1} +  ( \E |X(\tk)|^{-p_1} )^{p\b/p_1} \right ) \left ( \E |X(s) - X(\tk) |^{pp_1/(p_1 - p \b)} \right )^{(p_1 - p\b)/p_1} \no \\
 & \le C\D^{p/2}, \no
\end{align}
where Assumption \ref{H1} and Lemma \ref{Lem3.3_2} have been used.
Moreover,
\begin{align*}
M^{(d)}_k = \int_{\tk}^{\tkk}\Big ( f(X(s)) - f (X(\tk)) \Big )ds
        +\int_{\tk}^{\tkk}\Big ( f (X(\tk)) - \fd (X(\tk)) \Big )ds.
\end{align*}
By the \Holder inequality, \eqref{Eq94} and \eqref{Eq162_6}, we have
\begin{align}\label{Eq148_1}
&  \E \left [|M^{(d)}_k|^{p} \right ] \\
& \le C\D^{p-1} \int_{\tk}^{\tkk} \E \left  | f(X(s)) - f(X(\tk)) \right |^p ds
+  C\D^{p-1} \int_{\tk}^{\tkk} \E \left  | f(X(\tk)) - \fd (X(\tk)) \right |^p ds  \no \\
& \le  C\D^p  \Big (\frac{1}{\psi (\D )} +  \D^{p/2} \Big ). \no
\end{align}
By  the \BDG identity, we have
\begin{align} \label{Tem616}
& \E \left [|M^{(w)}_k|^p \right ]  \\
&  = \E \Big |  \int_{\tk}^{\tkk} \Big(    g(X(s)) - \gd (X(\tk)) \Big  ) dB(s) \Big |^p \no \\
& \le C \E \Big |  \int_{\tk}^{\tkk} \Big(    g(X(s)) - g (X(\tk)) \Big  ) dB(s) \Big |^p
+ C \E \Big |  \int_{\tk}^{\tkk} \Big(   g (X(\tk)) - \gd (X(\tk)) \Big  ) dB(s) \Big |^p \no \\
& \le C  \D^{p/2 -1} \int_{\tk}^{\tkk}  \E |    g(X(s)) - g(X(\tk))  |^p ds
+    C  \D^{p/2 -1}  \int_{\tk}^{\tkk}  \E |    g(X(\tk)) - \gd (X(\tk))   |^p ds . \no
\end{align}
Combining \eqref{Eq162_6} with \eqref{Eq94}, we conclude from \eqref{Tem616} that
\begin{align*}
\E \left [|M^{(w)}_k|^p \right ]&
\le
 C\D^{p/2}  \Big (\frac{1}{\psi (\D )} +  \D^{p/2} \Big ).
\end{align*}
The proof is finished.
$\Box$

We now consider the global discretization error between the true solution $X$ and the discretized process $X_\D$.
The following theorem reveals the optimal strong convergence rate of TEM in the sense of $L^p(\Omega; \R)$.
\begin{theorem}[Convergence order of TEM]\label{Theorem3.1}
Let Assumptions \ref{H0}, \ref{Assu1} and \ref{H1} hold with
\begin{align}\label{Cond3}
p_0 \ge 2(\a \ve \b) + 2 \a +2  \qu \textrm{and} \qu p_1 \ge  2(\a \ve \b) + 2 \b.
\end{align}
Let
\begin{align}\label{Cond5}
\displaystyle p \in \left [2,  \frac{p_0}{(2\a+1) \ve ( \a + \b +1)}
 \we \frac{p_1}{(\a+\b) \ve 2\b}  \right ].
\end{align}
Then, the truncated EM scheme \eqref{PTEM} by setting
\begin{align}\label{Set_gamma}
R(\D)= L_1 \D^{- \frac{1}{2 (\a \ve \b )}}
\end{align}
 has the property that
\begin{align}
\max_{0 \le k \le \lceil T/\D \rceil}\E[|X(\tk) - X_\D (\tk)|^{p}]  \le C \D^{p/2},\qu T>0, \label{Result1} \\
\max_{0 \le k \le \lceil T/\D \rceil}\E[|X(\tk) - \pd (X_\D (\tk))|^{p}]  \le  C \D^{p/2},  \qu T>0, \label{Result2}
\end{align}
where $C$ is a positive constant independent of $\D$.
\end{theorem}
\textbf{Proof.}
One notes that condition \eqref{Cond3} implies conditions \eqref{Cond1} and \eqref{Cond2}.
We observe from \eqref{Eq0} that
\begin{align}\label{Eq111}
X(\tkk) = X(\tk) + \fd (X(\tk))\D + \gd (X(\tk))\D B_k + M^{(d)}_k + M^{(w)}_k, \; k = 0,1,\cdots,
\end{align}
where $M^{(d)}_k$ and $M^{(w)}_k$ are defined in Lemma \ref{Lem3.5}.
Thus, we conclude from \eqref{PTEM2} and \eqref{Eq111} that
\begin{align}\label{Eq112}
X(\tkk) - X_\D(\tkk)  = &[X(\tk) - X_\D(\tk) ] + [\fd (X(\tk)) - \fd(X_\D(\tk))]\D  \no \\
&  + [\gd (X(\tk)) - \gd(X_\D(\tk))]\D B_k + M^{(d)}_k + M^{(w)}_k, \; k = 0,1,\cdots.
\end{align}
Denote by
\begin{align}\label{Tempp1}
A_k:= X(\tkk) - X_\D(\tkk) , \qu F_k := \fd (X(\tk)) - \fd(X_\D(\tk)), \qu G_k: = \gd (X(\tk)) - \gd(X_\D(\tk)).
\end{align}
Without loss of generality, we assume that $p \ge 2$ is an even integer and Condition \eqref{Cond5} is satisfied.
It follows from the binomial formula that
\begin{align}\label{Tempp2}
A_{k+1}^{p} & = (A_k + F_k \D + G_k \D B_k +  M^{(d)}_k + M^{(w)}_k )^{p}
 =  A_k^{p}+  \sum_{i=1}^{p/2}\Pi_{2i-1},
\end{align}
 where
 \begin{align}\label{Tempp3}
\Pi_i & = C_{p}^{i} A_k^{p- i}(F_k \D + G_k \D B_k +  M^{(d)}_k + M^{(w)}_k)^{i}, \no \\
& \qu + C_{p}^{i+1} A_k^{p- (i+1)}(F_k \D + G_k \D B_k +  M^{(d)}_k + M^{(w)}_k)^{i+1},  \qu i = 1,3,\cdots,  p/2-1.
 \end{align}
By the elementary inequality  \eqref{Inequality1},
 we have
\begin{align*}
\Pi_1 & = p A_k^{p-1} (F_k \D + G_k \D B_k +  M^{(d)}_k + M^{(w)}_k)   \no \\
& \qu + \frac{p(p-1)}{2} A_k^{p-2} (F_k \D + G_k \D B_k +  M^{(d)}_k + M^{(w)}_k)^2 \no \\
& \le \underbrace{\frac{p}{2} |A_k|^{p-2} (2 A_k  F_k \D + q_0 |G_k |^2 |\D B_k|^2 )}_{=:J_{1}}  +
\underbrace{C|A_k|^{p-1} | M^{(d)}_k | }_{=:J_{2}}  +
\underbrace{C|A_k|^{p-2} |M^{(w)}_k|^2}_{=:J_{3}}\no \\
& \qu + \underbrace{ C |A_k|^{p-2}|F_k|^2\D^2}_{=:J_{4}} +
 \underbrace{C |A_k|^{p-2}  |M^{(d)}_k|^2}_{=:J_{5}} +
 \underbrace{2 A_k^{p-1}G_k \D B_k + p A_k^{p-1}M^{(w)}_k}_{=:J_{6}} \no \\
& = \sum_{i=1}^{6}J_{i}.
\end{align*}
According to Lemma \ref{Lem1}, we have
\begin{align}\label{Temp621_2}
\E [J_1] & = \frac{p}{2} \E \Big [ |X(\tk) - X_\D(\tk)|^{p-2} \no \\
& \qu \K \Big (2[X(\tk) - X_\D(\tk)]  [\fd (X(\tk)) - \fd(X_\D(\tk))]\D  + q_0 [\gd (X(\tk)) - \gd(X_\D(\tk))]^2  |\D B_k|^2    \Big )  \Big ]\no \\
& \le C \D \E \left [ |X(\tk) - X_\D(\tk)|^{p} \right ].
\end{align}
By the elementary inequality \eqref{Inequality1} and Lemma \ref{Lem3.5},
we have
\begin{align}\label{Temp621_5}
\E [J_2] =  C \E \left [ |X(\tk) - X_\D(\tk)|^{p-1} |M^{(d)}_k | \right ]&
\le C  \D \E \left [ |X(\tk) - X_\D(\tk)|^{p}  \right ]
+ C \frac{\E \left [ |M^{(d)}_k|^{p}\right ]}{\D^{p-1}}  \no \\
& \le  C  \D \E \left [ |X(\tk) - X_\D(\tk)|^{p}  \right ] +  C\D  \Big (\frac{1}{\psi (\D )} +  \D ^{p/2}\Big ) .
\end{align}
Similarly,
\begin{align}\label{Temp621_6}
\E [J_3] =  C \E \left [ |X(\tk) - X_\D(\tk)|^{p-2} |M^{(w)}_k |^2 \right ]&
\le C  \D \E \left [ |X(\tk) - X_\D(\tk)|^{p}  \right ]
+ C \frac{\E \left [ |M^{(w)}_k|^{p}\right ]}{\D^{p/2-1}}  \no \\
& \le  C  \D \E \left [ |X(\tk) - X_\D(\tk)|^{p}  \right ] +  C\D  \Big (\frac{1}{\psi (\D )} +  \D ^{p/2}\Big ) ,
\end{align}
and
\begin{align}\label{Temp621_7}
\E [J_5] =  C \E \left [ |X(\tk) - X_\D(\tk)|^{p-2} |M^{(d)}_k |^2 \right ]&
\le C  \D \E \left [ |X(\tk) - X_\D(\tk)|^{p}  \right ]
+ C \frac{\E \left [ |M^{(d)}_k|^{p}\right ]}{\D^{p/2-1}}  \no \\
& \le  C  \D \E \left [ |X(\tk) - X_\D(\tk)|^{p}  \right ] +  C\D^{p/2+1}  \Big (\frac{1}{\psi (\D )} +  \D ^{p/2}\Big )  \no \\
& \le  C  \D \E \left [ |X(\tk) - X_\D(\tk)|^{p}  \right ] +  C\D^{}  \Big (\frac{1}{\psi (\D )} +  \D ^{p/2}\Big ) .
\end{align}
 Moreover, $\fd$ is globally Lipschitz continuous with Lipschitz constant $\varphi(\D)$, see \eqref{eq62}. Thus
\begin{align}\label{Temp621-9}
\E [J_4] & = C \D^2 \E \left [ |X(\tk) - X_\D(\tk)|^{p-2} |\fd (X(\tk)) - \fd(X_\D(\tk))|^2  \right ]  \no \\
& \le C \varphi (\D)^2 \D^2 \E \left [ |X(\tk) - X_\D(\tk)|^{p}  \right ]
\le C \D \E \left [ |X(\tk) - X_\D(\tk)|^{p}  \right ].
\end{align}
In addition, we have the following identity
\begin{align}\label{Temp621_1}
\E [J_6] & =  p \E \left [ |X(\tk) - X_\D(\tk)|^{p-1} \big [ \gd (X(\tk)) - \gd(X_\D(\tk))\big ] \D B_k \right ] \no \\
&
+ p \E \left [(X(\tk) - X_\D(\tk)|^{p-1}
\int_{\tk}^{\tkk}\Big ( g(X(s)) - \gd (X(\tk)) \Big )dB(s) \right ] =0.
\end{align}
From \eqref{Temp621_2}-\eqref{Temp621_1},
we have
 \begin{align}\label{Temp621-10}
\E [\Pi_1]
\le C \D \E \left [ |X(\tk) - X_\D(\tk)|^{p}  \right ]+ C\D^{}  \Big (\frac{1}{\psi (\D )} +  \D ^{p/2}\Big ).
\end{align}
 In the same fashion as \eqref{Temp621-10} is obtained, we also can show that
\begin{align}\label{Temp621_11}
 \sum_{i=2}^{p/2} \E  \left [ \Pi_{2i-1} \right ]\le C \D \E \left [ |X(\tk) - X_\D(\tk)|^{p}  \right ]+ C\D^{}  \Big (\frac{1}{\psi (\D )} +  \D ^{p/2}\Big ).
\end{align}
Thus, we conclude from \eqref{Tempp2} that
\begin{align}\label{Eq116_2}
&\E [|X(\tkk) - X_\D(\tkk)|^{p}] \\
&  \le (1 + C \D ) \E [|X(\tk) - X_\D(\tk)|^{p}]
  +  C \D \Big (\frac{1}{\psi (\D )} +  \D^{p/2} \Big ),\; k = 0,1,\cdots. \no
\end{align}
Using the discrete Gronwall inequality yields
\begin{align}\label{Eq149_1}
 \max_{0 \le k \le N} \E [|X(\tk) - X_\D(\tk)|^{p}]
& \le   C  \Big (\frac{1}{\psi (\D )} +  \D^{p/2} \Big ) \no \\
&  = C \Big (  \D^{\g (p_0 - p \a  - p)} + \D^{\g ( p_1 - p\b  +p)} + \D^{p/2} \Big ).
\end{align}
To balance the error terms, we set
  $\displaystyle \g =\frac{1}{ 2 ( \a  \ve  \b) }$. From this and \eqref{Cond5}, we have
 \begin{align}\label{Temp623}
\g (p_0 - p \a  - p) \ge p/2 \qu \textrm{and} \qu \g ( p_1 - p\b  +p) \ge p/2,
 \end{align}
  which implies that
\begin{align}\label{Eq149_10}
\sup_{0 \le t \le T} \E \left [ | f(X(t)) - \fd (X(t)|^{p} \ve | g(X(t)) - \gd (X(t)|^{p}
\ve | X(t) -  \pd (X(t))|^{p} \right ]  \le C \D^{p/2},
\end{align}
and  \eqref{Result1} holds.
Moreover,
by Lemma \ref{Lem1} and \eqref{Eq94_3} as well as \eqref{Result1}, we have
\begin{align}\label{Eq149_2}
& \E [|X(\tk) - \pd( X_\D(\tk)) |^{p}] \\
& \le C\E [|X(\tk) - \pd( X(\tk)) |^{p}]
+   C\E [|\pd( X(\tk)) - \pd( X_\D(\tk)) |^{p}]  \no \\
&  \le C\E [|X(\tk) - \pd( X(\tk)) |^{p}]
+   C\E [| X(\tk) -  X_\D(\tk) |^{p}] \no \\
& \le C  \D^{\g (p_0  - 2) } +   C  \D^{\g (p_1  + 2) } + C  \D^{p/2}  \le C \D^{p/2}, \; k = 0,1,\cdots, \no
\end{align}
which is the desired assertion \eqref{Result2}.
$\Box$

\begin{remark}\label{Remark1}
Let $ \g \in (0, \frac{1}{2(\a \ve \b)}]$.
 From  the proofs of Theorem \ref{Theorem3.1}, we note that if
 Assumptions \ref{H0}, \ref{Assu1} and \ref{H1} hold with
$ p_0 \ge  (4\a +2) $ and $  p_1 \ge 4 \b $,
%
then for  $p \in \left [2,  \frac{p_0}{2\a+1} \we \frac{p_1}{2\b} \right ]$,
\begin{align}\label{Eq626_2}
\max_{0 \le k \le \lceil T/\D \rceil}\E[|X(\tk) - \pd ( X_\D (\tk))|^{p}]
& \le  C \Big (  \D^{\g (p_0 - p \a  - p)} + \D^{\g ( p_1 - p\b  +p)} + \D^{p/2} \Big ).
\end{align}
 Comparing \eqref{Eq626_2} with \cite[Eq. 3.15]{Guo2018note}, we observe  that when
$1/x$ term varnishes in \eqref{Eq0}, then Theorem \ref{Theorem3.1} degenerates to  \cite[Theorem 3.6]{Guo2018note}.
\end{remark}

For application use, we show that our numerical solutions have some finite, or inverse, moments.
\begin{proposition}[Moment boundedness of TEM]\label{Cor3.1}
Let Assumptions \ref{H0}, \ref{Assu1} and \ref{H1} hold with
\begin{align}\label{Cond31}
p_0 \ge  4(\a \ve \b) +8 \qu \textrm{and} \qu p_1 \ge  4 (\a \ve \b) +2.
\end{align}
Let
\begin{align}\label{Cond7}
\displaystyle p \in \left [2,  \frac{p_0}{2(\a \ve \b)+4 }
 \we \frac{p_1}{2 (\a \ve \b) +1}  \right ]
\end{align}
and  $R(\D)= L_1 \D^{- \g}$ with
\begin{align}\label{Set_gamma1}
\g = \frac{1}{2 (\a \ve \b )+4}.
\end{align}
Then,
\begin{align}\label{Eq626_9}
\max_{0 \le k \le \lceil T/\D \rceil}\E[|X(\tk) - \pd ( X_\D (\tk))|^{p}]
& \le  C \D^{p/2}, \qu T>0.
\end{align}
Moreover,
\begin{enumerate}
  \item[(1).] For $ 2 \le q \le p(\a \ve \b) + 2p +1$,
   \begin{align}\label{Moment_Y_1_1}
\max_{0 \le k \le \lceil T/\D \rceil}\E[|\pd (X_\D (\tk))|^q] & <  \infty, \qu T>0;
\end{align}
  \item[(2).] For $ 1 \le q \le p(\a \ve \b) +1$,
   \begin{align}\label{Moment_Y_1_2}
\max_{0 \le k \le \lceil T/\D \rceil}\E[|\pd (X_\D (\tk))|^{-q}] & <  \infty, \qu T>0.
\end{align}
\end{enumerate}

\end{proposition}
\textbf{Proof.}
Let $p$ and $\g$ be the parameters satisfying \eqref{Cond7} and \eqref{Set_gamma1}, respectively.
Then,
\begin{align*}
\g (p_0 - p \a - p)& \ge \frac{p[ 2 (\a \ve \b) +4 - \a -1]}{ 2(\a \ve \b)+4}\ge \frac{[(\a \ve \b) +3]p}{2(\a \ve \b) +4}  \ge \frac{p}{2}, \\
 \g (p_1 - p \b +1) & \ge \frac{p[ 2 (\a \ve \b) +1 - \b +1]}{ 2(\a \ve \b)+4} \ge \frac{[(\a \ve \b) +2]p}{2(\a \ve \b) +4}=\frac{p}{2}.
\end{align*}
Therefore, \eqref{Eq626_9} follows from \eqref{Eq626_2}.
We now  prove (1).
For $ 2 \le q \le p$,  \eqref{Moment_Y_1_1}    follows from Theorem  \ref{Theorem3.1}, Assumption \ref{H1} and the following inequality
\begin{align*}
| \pd (X_\D (\tk)) |^p \le 2^{p-1} |X(\tk)|^p +  2^{p-1} |X(\tk)-\pd (X_\D (\tk))|^p.
\end{align*}
We now assume that
$p< q \le p_0$.
Write $R = R(\D)$ for simplicity.
Let
\begin{align*}
\pd (X_\D (\tk)) =  e^\D_k + X(\tk),
\end{align*}
where $e^\D_k: = \pd (X_\D (\tk)) - X(\tk)$.
Clearly, we have that
\begin{align}\label{Eq181}
|\pd (X_\D (\tk))|^q
= & |\pd (X_\D (\tk))|^q \II_{\{ |X(\tk)|>R \}}
 + |\pd (X_\D (\tk))|^q \II_{\{ |X(\tk)|\le R, \; |e^\D_k| \le 1 \}} \\
& + |\pd (X_\D (\tk))|^q \II_{\{ |X(\tk)|\le R, \; |e^\D_k| > 1 \}}. \no
\end{align}
Since, $ 1/ R \le \pd (X_\D (\tk)) \le R$, we conclude from Assumption \ref{H1} that
\begin{align}\label{Eq182}
& \E \Big  [|\pd (X_\D (\tk))|^q \II_{\{ |X(\tk)|>R \}} \Big  ] \no \\
& \le R^q \E \Big [\II_{\{ |X(\tk)|>R \}} \Big ] \le R^q \frac{\E \left [ |X(\tk)|^q \right ]}{R^q} \le \E \left [ |X(\tk)|^{p_0} \right ]   \le C.
\end{align}
By the elementary inequality and Assumption \ref{H1}, we have
\begin{align}\label{Eq183}
\E \left [ |\pd (X_\D (\tk))|^q \II_{\{ |X(\tk)| \le R, \; |e^\D_k| \le 1 \}} \right ]
& = \E \left [ | X(\tk) + e^\D_k|^q
\II_{\{ |X(\tk)| \le R, \; |e^\D_k| \le 1 \}}
\right ] \\
& \le C \Big ( \E [ | X(\tk) |^q ] +
\E [ | e^\D_k |^p ]
\II_{\{ |X(\tk)| \le R, \; |e^\D_k| \le 1 \}}
   \Big ) \no \\
 & \le C \Big ( \E [ | X(\tk) |^q ] + 1 \Big ) \le C. \no
\end{align}
We note that
\begin{align*}
|y|^q \le |x|^q + C_q (|x|^{q-1} + |y|^{q-1})|x-y|, \qu \f x, y \in \R.
\end{align*}
Therefore,
\begin{align}\label{Temp615}
& |\pd (X_\D (\tk))|^q \II_{\{ |X(\tk)| \le R, \; |e^\D_k| > 1 \}}  \\
& \le |X(\tk)|^q +
 C_q \Big ( |X(\tk)|^{q-1} +|\pd (X_\D (\tk))|^{q-1}  \Big )
 |e^\D_k|
 \II_{\{ |X(\tk)| \le R, \; |e^\D_k| > 1 \}}. \no \\
\end{align}
By the \Holder and  the Markov inequalities as well as Theorem \ref{Theorem3.1}, we have
\begin{align}\label{Temp1}
\E \left [ |e^\D_k|\II_{\{|e^\D_k| > 1 \}} \right ]
& \le \Big ( \E \Big  [ |e^\D_k  |^p \Big ] \Big )^{1/p}
\Big ( \E \left [ \II_{\{|e^\D_k| > 1 \}} \right ] \Big )^{1-1/p}   \\
& \le \Big ( \E \Big  [ |e^\D_k  |^p \Big ] \Big )^{1/p}
\Big ( \E \Big  [ |e^\D_k  |^p \Big ] \Big )^{1- 1/p}
=  \E \Big  [ |e^\D_k  |^p \Big ].\no
\end{align}
If we set
\begin{align}\label{Eq626_4}
 q \le \frac{p}{2\g}+1 = p [(\a \ve \b)+2]+ 1
\end{align}
%
such that  $ \g (q-1) \le p/2$, then we conclude from \eqref{Temp615} and \eqref{Temp1} that
\begin{align}\label{Eq184}
& \E \left [|\pd (X_\D (\tk))|^q
\II_{\{ |X(\tk)| \le R, \; |e^\D_k| > 1 \}} \right ] \\
& \le \E \Big [ |X(\tk)|^q \Big ]+ C R^{q-1} \E \left [ |e^\D_k|
\II_{\{ |X(\tk)| \le R, \; |e^\D_k| > 1 \}} \right ] \no \\
& \le C(1 + R^{q-1}) \E \left [ |e^\D_k|\II_{\{|e^\D_k| > 1 \}} \right]
\no \\
& \le  C R^{q-1} \E \left [ |e^\D_k|^p \right] \le C R^{q-1} \D^{p/2}= C \D^{p/2 - \g (q-1)} \le C. \no
\end{align}
From \eqref{Eq181}-\eqref{Eq184}, we conclude that \eqref{Moment_Y_1_1} holds.
\par
We now prove (2). For $q \ge 1$ We observe that
\begin{align}\label{Eq191}
& |\pd (X_\D (\tk))|^{-q} \\
= & |\pd (X_\D (\tk))|^{-q} \II_{\{ |X(\tk)|^{-1}>R \}}
 + |\pd (X_\D (\tk))|^{-q} \II_{\{ |X(\tk)|^{-1}\le R, \; |e^\D_k| \le R^{-2} \}}  \no \\
& + |\pd (X_\D (\tk))|^{-q} \II_{\{ |X(\tk)|^{-1}\le R, \; |e^\D_k| > R^{-2} \}}. \no
\end{align}
Since, $ 1/ R \le 1/ \pd (X_\D (\tk)) \le R$, we conclude from Assumption \ref{H1} that
\begin{align}\label{Eq192}
&\E \Big [|\pd (X_\D (\tk))|^{-q} \II_{\{ |X(\tk)|^{-1}>R \}} \Big ]
\le R^q \E \Big [\II_{\{ |X(\tk)|^{-1}>R \}} \Big ]  \\
& \le R^q \frac{\E \left [ |X(\tk)|^{-q} \right ]}{R^q} \le C_q. \no
\end{align}
By the elementary inequality and Assumption \ref{H1}, we have
\begin{align}\label{Eq193}
& \E \left [ |\pd (X_\D (\tk))|^{-q} \II_{\{ |X(\tk)|^{-1} \le R, \; |e^\D_k| \le R^{-2} \}} \right ] \\
& \le C \left ( \E \left [\Big |  \frac{1}{X(\tk)} \Big |^q\right ] +
\E \left [ \Big | \frac{1}{X(\tk)} - \frac{1}{\pd (X(\tk))} \Big |^q
\II_{\{ |X(\tk)|^{-1} \le R, \; |e^\D_k| \le R^{-2} \}} \right ]
 \right ) \no \\
 & \le C \left ( \E \left [\Big |  \frac{1}{X(\tk)} \Big |^q\right ] +
\E \left [ \frac{|e^\D_k|^q}{|X(\tk)|^q |\pd ( X_\D (\tk))|^q}
\II_{\{ |X(\tk)|^{-1} \le R, \; |e^\D_k| \le R^{-2} \}} \right ]
 \right ) \no \\
  & \le C \left ( \E \left [\Big |  \frac{1}{X(\tk)} \Big |^q\right ] +
 R^{2q}\E \left [ |e^\D_k|^q
\II_{\{ |X(\tk)|^{-1} \le R, \; |e^\D_k| \le R^{-2} \}} \right ]
 \right ) \no \\
& \le C \left ( \E \left [\Big |  \frac{1}{X(\tk)} \Big |^q\right ]
+ 1
\right ) \le C. \no
\end{align}
By the \Holder and Markov inequality, we have
\begin{align}\label{Temp11}
\E \left [ |e^\D_k|\II_{\{|e^\D_k| > R^{-2} \}} \right ]
& \le \Big ( \E \Big  [ |e^\D_k  |^p \Big ] \Big )^{1/p}
\Big ( \E \left [ \II_{\{|e^\D_k| > R^{-2} \}} \right ] \Big )^{1-1/p}  \no \\
& \le \Big ( \E \Big  [ |e^\D_k  |^p \Big ] \Big )^{1/p}
\Big ( R^{2p} \E \Big  [ |e^\D_k  |^p \Big ] \Big )^{1-1/p}
=  R^{2p-2} \E \Big  [ |e^\D_k  |^p \Big ].
\end{align}
 If we set
 \begin{align}\label{Eq627_1}
q \le \frac{p}{2\g}-2p+1 = p (\a \ve \b ) +1
           \end{align}
            such that   $\g (q + 2p -1) \le p/2$,
%
then
\begin{align}\label{Eq194}
& \E \left [|\pd (X_\D (\tk))|^{-q}
\II_{\{ |X(\tk)|^{-1} \le R, \; |e^\D_k| > R^{-2} \}} \right ] \\
& \le \E \Big [ |X(\tk)|^{-q} \Big ]+
 C \E \left [ \Big ( |X(\tk)|^{-q-1}  + |\pd (X(\tk))|^{-p-1}\Big )  |e^\D_k|
\II_{\{ |X(\tk)|^{-1} \le R, \; |e^\D_k| > R^{-2} \}} \right ]
  \no \\
& \le \E \Big [ |X(\tk)|^{-q} \Big ]+ C_q R^{q+1} \E \left [ |e^\D_k|
\II_{\{ |X(\tk)|^{-q} \le R, \; |e^\D_k| > R^{-2} \}} \right ] \no \\
& \le C(1 + R^{q+1}) \E \left [ |e^\D_k|\II_{\{|e^\D_k| > R^{-2} \}} \right]
\no \\
& \le  C (1 + R^{q+1})R^{2p-2}\E \left [ |e^\D_k|^{p} \right]  \no \\
& \le C R^{q+2p-1} \D^{p/2}
\le C \D^{p/2 - \g (q + 2p -1)}
\le C. \no
\end{align}
From \eqref{Eq191}-\eqref{Eq194}, we conclude that \eqref{Moment_Y_1_2} holds.
$\Box$ \par
From the above analysis, we observe that  in $L^2(\Omega; \R)$ sense if we set
$\g = \frac{1}{(p_0 - 2 \a -2) \we ( p_1 - 2 \b +2)}$, then
our convergence results  can be described in a more concise form, see the following corollary.
%
\begin{corollary}\label{Cor3}
Let Assumptions \ref{H0}, \ref{Assu1} and \ref{H1} hold with $p_0=p_1$, $\a=\b$
and $p_0 = ( 4\a +2)$.
Let $R(\D)= L_1 \D^{- \frac{1}{2 \a  }}$.
Then
\begin{align}
& \max_{0 \le k \le \lceil T/\D \rceil}\E[|X(\tk) - \pd (X_\D (\tk))|^{2}]  \le  C \D,  \qu T>0,  \label{Res628_1}
\end{align}
and
\begin{align}\label{Res628_2}
& \max_{0 \le k \le \lceil T/\D \rceil} \E\left [|\pd (X_\D (\tk))|^{\frac{p_0}{2} \ve 2} \right ] < \infty \qu T>0.
\end{align}
Moreover, if $ p_0  \ge 10$, then
\begin{align}
& \max_{0 \le k \le \lceil T/\D \rceil}   \E\left [|\pd (X_\D (\tk))|^{-\left ( \frac{p_0}{2} -4 \right )} \right ] < \infty, \qu T>0. \label{Res628_3}
\end{align}
%
\end{corollary}

\begin{remark}\label{Remark2}
It is worth mentioning how our results of  TEM  compares with that of  Zhan and Li  in \cite{Zhan2024acm},
where the authors  proved the  $L^2$-convergence of order $1/2$ of a truncated EM for super-linear SDE.
However, their method do not preserve the  property of positiveness
 that we are interested in
%
and their results do not reveal  the
boundedness of moment and inverse moment  of  the numerical solutions. This hinders the further application  of the truncated method for some important SDEs, such as Lamperti transformed CIR model, see Application \ref{cir_ap}.
\end{remark}
\subsection{Positivity-preserving Truncated  Milstein (TMil) scheme}
We now  come to the second numerical scheme, which is called truncated  Milstein (TMil) scheme. In order to show the
first-order strong rate of convergence for this scheme,
we need the following assumption, which is stronger than Assumption \ref{H0}.
%
\begin{assumption}\label{H2}
Let $f, g \in \mathcal C^2 (\R_+)$.
There are $\hat \a \ge 1$,
$\hat \b \ge  0 $, $\hat K_1 >0$, such that for any $x,y \in \R_+$ it holds
\begin{align}
|f''(x)|  & \le \hat K_1
\Big ( 1+ x^{\hat \a-1}+ x^{-\hat \b} \Big ), \label{Eq_H2} \\
|g''(x)| & \le \hat K_1
\Big ( 1+ x^{\hat \a/2 -1}+  x^{-\hat \b/2} \Big ).\label{Eq_H22}
\end{align}
\end{assumption}
For a possibly enlarged $C$ the following estimates are an immediate consequence of Assumption \ref{H2} and the mean value theorem: For any $x,y \in \R_+$ it holds
\begin{align}
&|f'(x)-f'(y)| \le C \Big ( 1+ x^{\hat \a -1}+ y^{\hat \a-1} + x^{- \hat \b} +  y^{- \hat \b}  \Big)|x-y|,  \label{Eq5_2} \\
&|f(x)-f(y)| \le C \Big ( 1+ x^{\hat \a}+ y^{\hat \a} + x^{- \hat \b} +  y^{- \hat \b}  \Big)|x-y|,  \label{Eq5_4} \\
&|g'(x) - g'(y)| \le C \Big ( 1+ x^{\hat \a/2 -1}+ y^{\hat \a/2-1} + x^{- \hat \b/2} +  y^{- \hat \b/2} \Big)|x-y|,  \label{Eq5_21} \\
&|g(x) - g(y)| \le C \Big ( 1+  x^{\hat \a/2}+ y^{\hat \a/2} + x^{- \hat \b/2} +  y^{- \hat \b/2} \Big)|x-y|,  \label{Eq5_41} 
\end{align}
and
\begin{align}
&|f'(x)|  \le C \Big ( 1+ x^{\hat \a}+ x^{ -\hat \b}  \Big) , \label{Eq5_5}  \\
&|f(x)|  \le C \Big ( 1+ x^{\hat \a+1}+ x^{ -\hat \b}  \Big),  \label{Eq5_3} \\
& |g'(x)| \le C \Big ( 1+ x^{\hat \a /2}+ x^{ -\hat \b/2} \Big)  \label{Eq5_51} , \\
& |g(x)| \le C \Big ( 1+ x^{\hat \a/2+1}+ x^{ -\hat \b/2} \Big),  \label{Eq5_31} 
\end{align}
As above, we verify under Assumption \ref{H2} that the mapping $g' \cdot g$ satisfies
the polynomial Lipschitz condition
\begin{align}\label{Eq5_7}
&|g' \cdot g (x) -g' \cdot g (y) | \le C \Big ( 1+ x^{\hat \a}+ y^{\hat \a} + x^{-\hat \b} + y^{-\hat \b}  \Big)|x-y|,
\end{align}
and
the polynomial growth bound
\begin{align}\label{Eq_5_6}
|g' \cdot g (x)| \le C
\Big ( 1+ x^{\hat \a+1}+ x^{-\hat \b} \Big ),
\end{align}
for any $x,y \in \R_+$.
We define the following notation for the stochastic increments:
\begin{align*}
I_{t,s}:= \int_{t}^{s} \int_{t}^{t_1} dB(t_2) d B(t_1) = \frac{1}{2} \Big ( [B(s) - B(t)]^2 - (s-t) \Big ), \qu 0 \le t <s.
\end{align*}
We now begin to construct a  positivity-preserving truncated Milstein method, which
 is defined by setting $\hat X_\D  (t_0)   = X_0$ and by the recursion
\begin{align}
  \hat Y_\D (\tk )  & =  \pd (\hat X_\D (\tk )),   \no  \\
 \hat X_\D (\tkk)  & = \hat X_\D (\tk ) + f (\hat Y_\D (\tk ) ) + g (\hat Y_\D (\tk ) ) \D B_k +\frac{1}{2}g' \cdot g (\hat Y_\D (\tk)) ( |\D B_k|^2 - \D  ) ,   \label{PMil}\\
  & = \hat X_\D (\tk ) + \fd (\hat X_\D (\tk ) ) + \gd (\hat X_\D (\tk ) ) \D B_k +
  \ggd(\hat X_\D (\tk))I_{\tk, \tkk},  \label{PMil2}
\end{align}
for $k = 0,1,2, \cdots$,  where  $\pd $ has been defined in \eqref{eq23} with a given parameter  $ \g \in \left (0, \frac{1}{2(\hat \a \ve \hat \b )} \right ]$,
 $\fd $ and $\gd$  are defined in  \eqref{eq2},
\begin{align}\label{Eq_gg}
\ggd(x):= g' \cdot g (\pd (x)), \qu \f x \in \R.
      \end{align}

The following lemma implies that the TMil method is stochastically $C$-stable  in the sense of   \cite[Definition 3.2]{BIK2016PEM}, and plays an important role in the convergence analysis of the TMil method.

\begin{lemma}\label{Lem4.3}
Consider mapping $\pd $  defined in \eqref{eq23} with   $ \g \in \left (0, \frac{1}{2(\hat \a \ve \hat \b )} \right ]$.
Let  Assumption \ref{Assu1} hold.
Then there exists a constant $K$ only depending on $q_0$ such that
\begin{align*}
&|(x-y) -\D (\fd (x) - \fd (y) ) |  +  q_0 \D  |\gd(x)-\gd(y)|^2 + q_0 \D^2  |\ggd(x)-\ggd(y)|^2  \\
&  \le (1 + C\D) |x-y|^2, \quad \f x,y \in \R.
\end{align*}
\end{lemma}
\textbf{Proof.} By \eqref{Eq5_7} and  \eqref{Eq2}, we have
\begin{align*}
|\ggd(x)-\ggd(y)|\le C( 1 + R^{ \hat \a} + R^{ \hat \b}) |\pd(x) - \pd (y)| \le C \D^{- \g (\hat \a \ve \hat \b)} |x-y|, \qu \f x,y \in \R.
\end{align*}
Thus, by Lemma \ref{Lem1} and \eqref{eq62}, we have
\begin{align*}
&|(x-y) -\D (\fd (x) - \fd (y) ) |  +  q_0 \D  |\gd(x)-\gd(y)|^2 + q_0 \D^2  |\ggd(x)-\ggd(y)|^2  \\
& = |x-y|^2  + 2 \D \lan x-y,  \fd (x) - \fd (y) \ran  +  q_0 \D  |\gd(x)-\gd(y)|^2 \\
&  \qu + \D^2 |\fd (x) - \fd (y)|^2 + q_0 \D^2  |\ggd(x)-\ggd(y)|^2 \\
& \le  (1 + C \D) |x - y|^2 + C \D ^{2 - 2\g (\hat \a \ve \hat \b)} |x-y|^2  \\
& \le (1 + C \D) |x - y|^2.
\end{align*}
Thus, the proof is finished. $\Box$
\par
The following lemma is a consequence of the polynomial growth bound \eqref{Eq5_3}-\eqref{Eq5_51}.
\begin{lemma}\label{Lem712}
Let Assumptions \ref{H2}, \ref{Assu1} and \ref{H1} hold. Then,
\begin{align}
& \sup_{0 \le t \le T} \E \left [| \mathbb L f(X(t)) |^2 \right ]
\ve \sup_{0 \le t \le T} \E \left [| \mathbb L g(X(t)) |^2 \right ]
 \le C \Big (  1+ \sup_{0 \le t \le T} \E \big [ |X(t)|^{4\hat \a +2}\big ] +
  \sup_{0 \le t \le T} \E\big [ |X(t)|^{- 4\hat \b}\big ]
\Big ) ,   \label{Tem712-2}  \\
 & \sup_{0 \le t \le T} \E \left [ | f'(X(t))  g(X(t)) |^2 \right ]
  \le C \Big (  1+ \sup_{0 \le t \le T} \E \big [ |X(t)|^{3\hat \a +2}\big ] +
  \sup_{0 \le t \le T} \E\big [ |X(t)|^{- 3\hat \b}\big ]
\Big ) ,   \label{Tem712-3}
\end{align}

\end{lemma}
From \eqref{Eq148_1}, we have
\begin{align}\label{Tem_712-1}
\E \left [ \Big | \int_{t}^{t+\D}\Big ( f(X(s)) - f (X(t)) \Big )ds \Big |^2 \right ]
\le C \Big (  1+ \sup_{0 \le t \le T} \E \big [ |X(t)|^{4\hat \a +2}\big ] +
  \sup_{0 \le t \le T} \E\big [ |X(t)|^{- 4\hat \b}\big ]
\Big ) \D^3.
\end{align}
However, if we insert the conditional expectation with respect to the $\sigma$-field $\mathcal F_t$, the order convergence indicated by \eqref{Tem_712-1} can be increased, see \eqref{Eq713_1}.

\begin{lemma}\label{lem3.55}
Let Assumptions \ref{H2}, \ref{Assu1} and \ref{H1} hold with
\begin{align*}
p_0 \ge 4 \hat \a +2  \qu \textrm{and} \qu p_1 \ge  4\hat \b.
\end{align*}
Then for any $ t \in [0,T]$,
\begin{align}
&\E \left [ \Big  | \E_t \Big [ \int_{t}^{t +\D} \Big (  f(X(s)) - f(X(t)) \Big ) ds   \Big ] \Big   |^2 \right ]
  =\E \left [ \Big  |
\int_{t}^{t +\D} \Big (\int_{t}^{s}  \mathbb L f(X(u)) du \Big ) ds
 \Big |^2 \right ]
\le C \D^4,  \label{Eq713_1} \\
& \E \left [ \Big  | \int_{t}^{t+\D} \Big ( \int_{t}^{s} f'(X(u))  g(X(u))dB(u) \Big ) ds \Big |^2 \right ]  \le C \D^3,  \label{Eq713_2} \\
& \E \left [ \Big  | \int_{t}^{t+\D} \Big (
g(X(s)) - g(X(t)) \Big ) dB(s) - g'(X(t))\cdot g(X(t)) I_{t, t + \D} \Big |^2 \right ]  \le C \D^3. \label{Eq713_3}
\end{align}
\end{lemma}
\textbf{Proof.}
By the \Ito formula, we have
\begin{align}
f(X(s))- f(X(t)) = \int_{t}^{s}  \mathbb L f(X(u))du+
\int_{t}^{s} f'(X(u))  g(X(u))dB(u), \qu 0 \le t \le s.  \label{Tmp712-5}
\end{align}
Thus,
\begin{align}\label{Tmp712-7}
\E_t \Big [ \int_{t}^{t +\D} \Big (  f(X(s)) - f(X(t)) \Big ) ds   \Big ]
= \int_{t}^{t +\D} \Big (\int_{t}^{s}  \mathbb L f(X(u)) du \Big ) ds
\end{align}
By \eqref{Tem712-2},  we have
\begin{align*}
\E \left [ \Big | \int_{t}^{t +\D} \Big (\int_{t}^{s}  \mathbb L f(X(u)) du \Big ) ds \Big |^2 \right ]
\le \D^2  \int_{t}^{t +\D} \Big (\int_{t}^{s} \E[ |\mathbb L f(X(u))|^2] du \Big ) ds
\le C \D^4,
\end{align*}
which implies that \eqref{Eq713_1} holds.
\par
Moreover,
by the \Ito isometry and \eqref{Tem712-3}, we also have
\begin{align*}
& \E \left [ \Big  | \int_{t}^{t+\D} \Big ( \int_{t}^{s} f'(X(u))  g(X(u))dB(u) \Big ) ds \Big |^2 \right ]
 \\
& \le \D \int_{t}^{t+\D}\Big ( \E \left [ \Big | \int_{t}^{s} f'(X(u))  g(X(u))dB(u) \Big |^2  \right ] \Big )  ds \\
& = \D \int_{t}^{t+\D}  \int_{t}^{s} \E \left [ | f'(X(u))  g(X(u)) |^2 \right ] du   ds
\le C \D^3,
\end{align*}
which is the desired assertion  \eqref{Eq713_2}.

It remains to show \eqref{Eq713_3}.
Again, by the \Ito isometry, we have
\begin{align}\label{Tmp713_10}
& \E \left [ \Big  | \int_{t}^{t+\D} \Big (
g(X(s)) - g(X(t)) \Big ) dB(s) - g'\cdot g(X(t)) I_{t, t + \D} \Big |^2 \right ]  \no \\
& = \E \left [ \Big  | \int_{t}^{t+\D} \Big (
g(X(s)) - g(X(t))  - g'\cdot g(X(t)) [B(s) - B(t)] \Big ) dB(s)\Big |^2 \right ]  \no \\
& =  \int_{t}^{t+\D}  \Gamma(s) ds,
\end{align}
where
$$\Gamma(s):=\E \left [  \Big  |
g(X(s)) - g(X(t))  - g'\cdot g(X(t)) [B(s) - B(t)] \Big |^2 \right ] . $$
Thus, the assertion \eqref{Eq713_3} is proved if there is a constant $C$ independent of $s,t$ and $\D$ such that
\begin{align}\label{Eq713_20}
\Gamma(s) \le C \D^2, \qu \f s \in [t, t + \D].
\end{align}
Again, by the \Ito formula, we get
\begin{align*}
g(X(s)) - g(X(t))= \int_{t}^{s } \mathbb L g(X(u))du + \int_{t}^{s } g' \cdot g(X(u)) dB(u), \qu 0 \le t \le s.
\end{align*}
Thus,
\begin{align}\label{eq7_2}
\Gamma(s) & = \E \left [ \Big | \int_{t}^{s}  \mathbb L g(X(u))du +
 \int_{t}^{s} \Big ( g'\cdot g(X(u)) - g'\cdot g(X(t))  \Big )   dB(s)
   \Big |^2 \right ]  \no \\
   & \le
  \underbrace{   2 \E \left [ \Big | \int_{t}^{s}  \mathbb L g(X(u))du \Big | ^2 \right ]}
  _{:= \Gamma_1(s)}
   +  \underbrace{  2 \E \left [ \Big | \int_{t}^{s} \Big ( g'\cdot g(X(u)) - g'\cdot g(X(t))  \Big ) dB(s)  \Big | ^2 \right ]}_{:= \Gamma_2(s)}.
\end{align}
From \eqref{Tem712-2}, we have
\begin{align}\label{eq7_4}
\Gamma_1(s) \le 2 \D \int_{t}^{s}  \E \left [ | \mathbb L g (X(u))|^2 \right ] ds \le C \D^2.
\end{align}
By the \Ito isometry, we get
\begin{align}\label{eq7_6}
\Gamma_2(s)=  2 \int_{t}^{s} \E \left [ \Big |  \Big ( g'\cdot g(X(u)) - g' \cdot g(X(t))  \Big )\Big | ^2 \right ] ds .
\end{align}
In the similar fashion as \eqref{Eq162_6} is obtained, we also can show that
\begin{align}\label{eq7_5}
& \E \left [ \Big |  \Big ( g'\cdot g(X(u)) - g'\cdot g(X(t))  \Big )\Big | ^2 \right ]
 \no \\
& \le C \Big ( 1  +
\sup_{0 \le t \le T} \E \big [ |X(t)|^{4\hat \a +2}\big ] +
\sup_{0 \le t \le T} \E\big [ |X(t)|^{- 4\hat \b}\big ]
\Big ) \D , \qu  t \le u \le s \le t +\D.
\end{align}
Thus, inserting \eqref{eq7_6} and \eqref{eq7_4} into \eqref{eq7_2}, we conclude that
$\Gamma(s) \le C \D^2, \qu \f s \in [t, t + \D]$, which completes the proof of \eqref{Eq713_3}.
 $\Box$

\begin{lemma}\label{Lem4.5}
Let Assumptions \ref{H2}, \ref{Assu1} and \ref{H1} hold with
\begin{align}\label{Cond10}
p_0 \ge  4 (\hat \a \ve \hat \b ) + 2 \hat \a +2 \qu \textrm{and} \qu  p_1 \ge  4 (\hat \a \ve \hat \b ) + 2 \hat \b -2 .
\end{align}
Then
\begin{align}
\E \left [ \big |\E_{\tk} [\hat M^{(d)}_k] \big |^2 \right ]  & \le C\D^4, \quad  k = 0,1, \cdots,\label{Eq148_11} \\
\E \left [ \big |\hat M^{(w)}_k \big |^2 \right ] & \le C\D^3, \quad k = 0,1, \cdots, \label{Eq148_22}
\end{align}
where
\begin{align*}
& \hat M^{(d)}_k: = \int_{\tk}^{\tkk}\Big ( f(X(s)) - \fd ( X(\tk)) \Big )ds, \\
& \hat M^{(w)}_k: = \int_{\tk}^{\tkk}\Big ( g(X(s)) - \gd ( X(\tk)) \Big )dB(s) - \ggd ( X(\tk)) I_{\tk, \tkk},
\end{align*}
with $R(\D)= L_1 \D^{- \frac{1}{2 (\hat \a \ve \hat \b )}}$ and $I_{\tk, \tkk} = \frac{1}{2} (|\D B_k|^2 - \D)$.
\end{lemma}
\textbf{Proof.}
Let condition \eqref{Cond10} hold and set $\g = \frac{1}{2 (\hat \a \ve \hat \b)}$. Then
\begin{align}
\g ( p_0 - 2 \hat \a -2 ) \ge 2 \qu \textrm{and} \qu  \g ( p_1 - 2 \hat \b  +2) \ge 2.\label{tmp713_21}
\end{align}
Thus, we conclude from Lemma \ref{Lem3.1} that
\begin{align}\label{Tem_713_11}
& \sup_{0 \le t \le T} \E \left [ | f(X(t)) - \fd (X(t)|^p \right ] \ve
\sup_{0 \le t \le T} \E \left [ | g(X(t)) - \gd (X(t)|^p \right ]  \no \\
& \le
C \Big (  \D^{ \g(p_0 - 2\hat \a -2)} + \D^{ \g(p_1 - 2\hat \b +2)} \Big ) \le C \D^2.
\end{align}
Moreover, we have the following decomposition
\begin{align*}
\hat M^{(d)}_k = \int_{\tk}^{\tkk}\Big ( f(X(s)) - f (  X(\tk)) \Big )ds
+  \Big ( f ( X(\tk)) - \fd ( X(\tk)) \Big )\D.
\end{align*}
Therefore,
\begin{align*}
&\E \left [ \big |\E_{\tk} [\hat M^{(d)}_k] \big |^2 \right ] \\
& \le
2 \E \left [ \Big  | \E_{\tk} \Big [ \int_{\tk}^{\tk +\D} \Big (  f(X(s)) - f(X(\tk)) \Big ) ds   \Big ] \Big   |^2 \right ]  +  2 \D ^2 \E \left [ \Big  | f ( X(\tk)) - \fd ( X(\tk))  \Big |^2 \right ]  \\
& \le C\D^4 + C \D^{4 }  = C\D^4,
\end{align*}
where \eqref{Eq713_1} and \eqref{Tem_713_11} have  been used.
In addition, \Ito isometry implies that
\begin{align}\label{tmp222}
\E \left [ \big |\hat M^{(w)}_k \big |^2 \right ]
& =\E \left [ \Big  | \int_{\tk}^{\tkk}  \Big (  g(X(s)) - \gd ( X(\tk))-
\ggd ( X(\tk)) [B(s)- B(\tk)] \Big ) dB(s)   \Big |^2 \right ] \no  \\
& = \int_{\tk}^{\tkk} \hat \Gamma (s) ds,
\end{align}
where
\begin{align*}
\hat \Gamma (s):= \E  \left [ \big | g(X(s)) - \gd ( X(\tk))-
\ggd ( X(\tk)) [B(s)- B(\tk)] \big |^2 \right ] .
\end{align*}
By  \eqref{Eq713_20} and \eqref{Eq94}, we have
\begin{align*}
\hat \Gamma (s) &\le 3 \E \left [ \big | g(X(s)) - g(X(\tk)) - g' \cdot g (X(\tk)) [B(s)- B(\tk)]  \big |^2  \right ] \\
& \qu + 3 \E \left [ \big | g(X(\tk)) - \gd (X(\tk)) \big|^2 \right ]
+ 3   \E \left [ \big | g' \cdot g (X(\tk))  - \ggd (X(\tk)) \big|^2 \big | B(s) - B(\tk) \big |^2\right ] \\
& \le C\D^2 + C \Big (  \D^{\g(p_0 - 2\hat \a -2)} + \D^{\g(p_1 - 2\hat \b +2)} \Big ) + C \D \E \left [ \big | g' \cdot g (X(\tk))  - \ggd (X(\tk)) \big|^2 \right ], \\
& \le C\D^2 + C \D \E \left [ \big | g' \cdot g (X(\tk))  - \ggd (X(\tk)) \big|^2 \right ],
\end{align*}
where the last step follows from  \eqref{tmp713_21}.
In a similar fashion as \eqref{Eq94} was proved, we also can show that
\begin{align*}
\sup_{ 0 \le t \le T} \E \left [ \big | g' \cdot g (X(t))  - \ggd (X(t)) \big|^2 \right ]
\le C \Big (  \D^{\g(p_0 - 2\hat \a -2)} + \D^{\g(p_1 - 2\hat \b +2)} \Big ) \le C\D^2,
\end{align*}
where \eqref{tmp713_21} has been used.
Thus, $$ \hat \Gamma(s) \le C \D^2. $$
Inserting this into \eqref{tmp222} completes the proof of \eqref{Eq148_22}.
 $\Box$
 \par
 The following   theorem shows that the TMil achieves the optimal mean-square convergence order $1$.
\begin{theorem}[Convergence order of TMil]\label{Theorem4.1}
Let Assumptions \ref{H2}, \ref{Assu1} and \ref{H1} hold with
\begin{align}\label{Cond33}
p_0 \ge  4 (\hat \a \ve \hat \b ) + 2 \hat \a +2 \qu \textrm{and} \qu  p_1 \ge  4 (\hat \a \ve \hat \b ) + 2 \hat \b  .
\end{align}
Then, the TMil scheme \eqref{PMil} by setting
\begin{align}\label{Set_gamma4}
R(\D)= L_1 \D^{- \frac{1}{2(\hat \a \ve \hat \b )}}
\end{align}
 has the property that
\begin{align}
\max_{0 \le k \le \lceil T/\D \rceil}\E[|X(\tk) - \hat X_\D (\tk)|^{2}]  \le C \D^{2},\qu T>0, \label{Result11} \\
\max_{0 \le k \le \lceil T/\D \rceil}\E[|X(\tk) - \pd (\hat X_\D (\tk))|^{2}]  \le  C \D^{2},  \qu T>0, \label{Result22}
\end{align}
where $C$ is a positive constant independent of $\D$.
\end{theorem}
\textbf{Proof.}
We observe from \eqref{Eq0} that
\begin{align}\label{Eq111_2}
X(\tkk) = X(\tk) + \fd (X(\tk))\D + \gd (X(\tk))\D B_k + \ggd ( X(\tk)) I_{\tk, \tkk} +\hat M^{(d)}_k + \hat M^{(w)}_k,
\end{align}
where $\hat M^{(d)}_k$ and $\hat M^{(w)}_k$ are defined in Lemma \ref{Lem4.5}.
Thus, we conclude from \eqref{PMil2} and \eqref{Eq111_2} that
\begin{align}\label{Eq112_2}
& \hat e^\D_{k+1}:= X(\tkk) - \hat X_\D(\tkk) \\
&  = [X(\tk) - \hat X_\D(\tk) ] + [\fd (X(\tk)) - \fd(\hat X_\D(\tk))]\D  + [\gd (X(\tk)) - \gd(\hat X_\D(\tk))]\D B_k\no \no \\
&  \qu +  \big [ \ggd ( X(\tk)) - \ggd (\hat  X_\D(\tk)) \big ]I_{\tk, \tkk}    + \hat M^{(d)}_k + \hat M^{(w)}_k.
\end{align}
By the orthogonality of the conditional expectation it holds
\begin{align*}
\E \left [\big |e^\D_{k+1}\big |^2 \right ] = \E \left [\big |\E_{\tk} [e^\D_{k+1}] \big |^2 \right ]
+ \E \left [\big | e^\D_{k+1} - \E_{\tk} [e^\D_{k+1}] \big |^2 \right ]
\end{align*}
Thus, by the elementary inequality
\begin{align}\label{Inequaltiy2}
(a+b)^2 \le (1 + \eps)a^2 + (1+1/\eps)b^2, \qu a, b, \eps >0,
\end{align}
and  $$\E \big [ |I_{\tk,\tkk}|^2\big ] = \frac{1}{2} \E \Big[|\D B_k - \D|^2 \Big ] = \frac{\D^2}{2},$$
as well as Lemma \ref{Lem4.3},
  we  have
\begin{align} \label{E_713}
& \E \left [ \big  | X(\tkk) - \hat X_\D(\tkk)  \big |^2 \right ]  \\
& = \E \left [ \big  |  [X(\tk) - \hat X_\D(\tk) ] +  [\fd (X(\tk)) - \fd(\hat X_\D(\tk))]  \D
+ \E_{\tk}[\hat M^{(d)}_k] \big |^2 \right ] +   \no \\
& \qu + \E \Big [ \big  | [\gd (X(\tk)) - \gd(\hat X_\D(\tk))]\D B_k   +
 \big [ \ggd ( X(\tk)) - \ggd (\hat  X_\D(\tk)) \big ]I_{\tk,\tkk} \no\\
& \qu + \hat M^{(w)}_k  +  \big ( \hat M^{(d)}_k -  \E_{\tk}[\hat M^{(d)}_k] \big ) \big |^2
\Big ]   \no\\
& \le (1 + \D)  \E \left [ \big  |  [X(\tk) - \hat X_\D(\tk) ] +  \D [\fd (X(\tk)) - \fd(\hat X_\D(\tk))]  \big |^2 \right ] \no\\
& \qu + q_0 \D \E \big [ |\gd (X(\tk)) - \gd(\hat X_\D(\tk))|^2 \big ] +
q_0 \D^2 \E \big [ |\ggd (X(\tk)) - \ggd (\hat X_\D(\tk))|^2 \big ] \no\\
& \qu+ \left (1 + \frac{1}{\D} \right )  \E \left [ \big  | \E_{\tk}[\hat M^{(d)}_k] \big |^2 \right ]
+ C \E \big [ \big |\hat M^{(w)}_k  \big |^2 \big ]
+
C \E \big [ \big | \hat M^{(d)}_k -  \E_{\tk}[\hat M^{(d)}_k] \big |^2 \big ] \no \\
& \le (1 + C \D) \E \left [ \big  | X(\tk) - \hat X_\D(\tk)  \big |^2 \right ] + \Xi_k,
\end{align}
where
 \begin{align*}
 \Xi_k:=\left (1 + \frac{1}{\D} \right ) \E \left [ \big  | \E_{\tk}[\hat M^{(d)}_k] \big |^2 \right ]
+ C \E \big [ \big |\hat M^{(w)}_k  \big |^2 \big ]
+C \E \big [ \big | \hat M^{(d)}_k -  \E_{\tk}[\hat M^{(d)}_k] \big |^2 \big ].
    \end{align*}
According to \eqref{Eq713_2}, we have
\begin{align*}
\E \big [ \big | \hat M^{(d)}_k -  \E_{\tk}[\hat M^{(d)}_k] \big |^2 \big ]= \E \left [ \Big  | \int_{\tk}^{\tk+\D} \Big ( \int_{\tk}^{s} f'(X(u))  g(X(u))dB(u) \Big ) ds \Big |^2 \right ]  \le C \D^3.
\end{align*}
Combining this and Lemma \ref{Lem4.5}, we get
$$ \Xi_k \le C \D^3.$$
Inserting this into \eqref{E_713}, we have
\begin{align*}
\E \left [ \big  | X(\tkk) - \hat X_\D(\tkk)  \big |^2 \right ]
\le (1 + C \D) \E \left [ \big  | X(\tk) - \hat X_\D(\tk)  \big |^2 \right ]  + C\D^3, \qu k = 0,1, \cdots.
\end{align*}
Using the discrete Gronwall inequality yields the desired assertion \eqref{Result11}.
It remains to prove \eqref{Result22}.
From Lemma \ref{Lem3.1}, we get immediately that
\begin{align*}
\max_{0 \le k \le \lceil T/\D \rceil}\E \left  [\big |X(\tk) - \pd ( X (\tk)) \big |^{2}\right ]
\le C \big ( \D^{\g (p_0 - 2)} + \D^{\g (p_1 + 2)} \big ) \le C \D^2.
\end{align*}
Combining this and Lemma \ref{Lem1} as well as the triangle inequality, we obtain the assertion \eqref{Result22}. Thus, the proof is finished.
 $\Box$

If we apply  Proposition \ref{Cor3.1}  with $p=2$, $\g = \frac{1}{2 (\hat \a \ve \hat \b)}$ and $\E \left [ |e^\D_k|^{p} \right]  \le C \D^2$, then we have
the following moment boundedness of the TMil solutions.
 \begin{proposition}[Moment boundedness of TMil]\label{Cor3.10}
Under the same assumption as Theorem  \ref{Theorem4.1},
the TMil scheme \eqref{PMil}
by setting \eqref{Set_gamma4} has the property that
\begin{enumerate}
  \item[(1).] For $ 2 \le q \le 4(\hat \a \ve \hat \b)  +1
  \le (p_0 - 2 \hat \a -1) \we (p_1 - 2 \hat \b -1)$,
   \begin{align}\label{Moment_Y_1_11}
\max_{0 \le k \le \lceil T/\D \rceil}\E[|\pd (\hat X_\D (\tk))|^q] & <  \infty, \qu T>0;
\end{align}
  \item[(2).] For $ 1 \le q \le 4(\hat \a \ve \hat \b) -3
  \le (p_0 - 2 \hat \a -5) \we (p_1 - 2 \hat \b -6)
  $,
   \begin{align}\label{Moment_Y_1_21}
\max_{0 \le k \le \lceil T/\D \rceil}\E[|\pd (\hat X_\D (\tk))|^{-q}] & <  \infty, \qu T>0.
\end{align}
\end{enumerate}

\end{proposition}

\section{Applications}\label{Sec4}
We now apply our results to some SDEs in mathematical finance.
\subsection{3/2 model}\label{32model}
The 3/2 process  $X$ is the solution to \eqref{model1}
and is strictly positive almost surely. Introduce the quantity
\begin{align}\label{Tem616_3}
\lambda:= 2 + \frac{2c_1}{\sigma^2}.
\end{align}
Existence and uniqueness can be retrieved from the properties of the Feller diffusion, and
\begin{align}\label{Temp616_4}
\sup_{0 \le t \le T} \E |X(t)|^p < \infty, \qu \f p < \lambda,
\end{align}
see \cite[p. 1009]{FCJ2016SIAM}. Clearly,  Assumption \ref{H0} is satisfied with exponents $\a=1$ and $\b=0$.  If condition
\begin{align}\label{Cond4}
\lambda \ge \frac{9}{8}q_0 +2, \qu \textrm{i.e., }\qu \frac{4c_1}{\sigma^2} \ge \frac{9}{4}q_0
\end{align}
holds,
then
\begin{align*}
2f'(x) + q_0 |g'(x)|^2 = 2 c_1 c_2 - \left (4c_1 - \frac{9}{4}q_0\sigma^2 \right )x^2 \le 2c_1c_2,
\end{align*}
which means that Assumption \ref{Assu1} holds.
If $\lambda >6$, we choose $p_0 \in [6, \l)$, $q_0=2$ and fix $\displaystyle \g = \frac{1}{p_0 - 4 }$, so that conditions \ref{Cond3} and \ref{Cond4} are satisfied.
The following results reveal that the TEM for 3/2 model has $L^2$-order of $1/2$ ,
$\E \left [ |\pd (X_\D (\tk ))|^p \right ]$  and
$\E \left [ |\pd (X_\D (\tk ))|^{-p} \right ]$ are  finite from Proposition \ref{Cor3}.
\begin{corollary}[Convergence rate and moment boundedness of TEM for 3/2 model] \label{cor31}
Let $\lambda >6$, i.e.,
$\sigma^2 < c_1/2$.
Then
the truncated EM scheme \eqref{PTEM} for model \eqref{model1}
by setting
$$ R(\D)= L_1 \D^{- \frac{1}{\l - 4}}$$ has the property that
\begin{align}\label{Tem61610}
\max_{0 \le k \le \lceil T/\D \rceil}\E[|X(\tk) - \pd (X_\D (\tk))|^2]  \le  C \D.
\end{align}
\begin{enumerate}
  \item[(1).] For any  $p \in [2, \l -3)$,
   \begin{align}\label{Moment_32}
\max_{0 \le k \le \lceil T/\D \rceil}\E[|\pd (X_\D (\tk))|^p] & <  \infty.
\end{align}
  \item[(2).] If $\l >8$, then for any  $p \in [1, \l -7)$,
   \begin{align}\label{Moment_321}
\max_{0 \le k \le \lceil T/\D \rceil}\E[|\pd (X_\D (\tk))|^{-p}] & <  \infty .
\end{align}
\end{enumerate}
\end{corollary}

Note that for model \eqref{model1},  Assumption \ref{H2} is satisfied with  $\hat \a =1$ and $\hat \b =0$.  The following assertion follows  directly from an application of
Theorem \ref{Theorem4.1} and Proposition \ref{Cor3.10}.
%
\begin{corollary}[Convergence rate of TMil for 3/2 model] \label{cor32}
Let $\lambda >8$, i.e.,
$\sigma^2 < c_1/3$.
Then
the truncated  Milstein scheme \eqref{PMil} for model \eqref{model1}
by setting
$$ R(\D)= L_1 \D^{- \frac{1}{2}}$$ has the property that
\begin{align}
\max_{0 \le k \le \lceil T/\D \rceil}\E[|X(\tk) - \hat X_\D (\tk)|^2]  \le  C \D^2, \label{Tem61610} \\
\max_{0 \le k \le \lceil T/\D \rceil}\E[|X(\tk) - \pd (\hat X_\D (\tk))|^2]  \le  C \D^2.\label{Tem61611}
\end{align}
Moreover,
\begin{align*}
\max_{0 \le k \le \lceil T/\D \rceil}\E[|\pd (\hat X_\D (\tk))|^{\l -3 }]
\ve \max_{0 \le k \le \lceil T/\D \rceil}\E[|\pd (\hat X_\D (\tk))|^{-{(\l -7)}}]
 <  \infty.
\end{align*}
\end{corollary}

\subsection{A\"{i}t-Sahalia model} \label{AIT_app}
Let
$X$ be  the solution to A\"{i}t-Sahalia interest rate model \eqref{Ait1}.
If $\kappa +1 > 2 \theta$, then
there exists a strong solution on $(0, \infty)$;  $\sup_{0 \le t \le T} \E |X(t)|^p$ and $\sup_{0 \le t \le T} \E |X(t)|^{-p}$ are finite for any $p>0$, see \cite[Lemma 2.1]{Mao2011NA_Sahalia}. In other words, Assumption \ref{H1} holds.
Moreover, Assumptions \ref{H0} and \ref{Assu1} also hold for $ \a= \kappa-1$, $\b =2$ and $q_0>2$, see \cite[p. 11]{Deng2023BIT}.
The following results
follows from
 Proposition \ref{Cor3.1}.
\begin{corollary}[Convergence rate of TEM for A\"{i}t-Sahalia model]\label{Cor33}
Let $ \k +1 > 2 \theta$ with $\k, \theta >1$.
Then, for any $p \ge 2$,
the TEM scheme \eqref{PTEM} for model \eqref{Ait1}
by setting
$$ R(\D)= L_1 \D^{- \frac{1}{ (2 \k  +2 )\ve 8}}$$ has the property that
\begin{align}
\max_{0 \le k \le \lceil T/\D \rceil}\E \left [|X(\tk) - \pd (X_\D (\tk))|^{p}\right ] & \le  C \D^{p/2}, \qu T >0, \label{Tem61610}\\
\max_{0 \le k \le \lceil T/\D \rceil}\E\Big [|\pd (X_\D (\tk))|^{p[(\k +1) \ve 4]+1} \Big ] & <  \infty, \qu T >0,  \label{tmp617_1} \\
\max_{0 \le k \le \lceil T/\D \rceil}
\E\left [ \frac{1}{ |\pd (X_\D (\tk))|^{p[(\k-1) \ve 2] +1}}   \right ]
& <  \infty, \qu T >0.  \label{tmp617_2}
\end{align}
\end{corollary}
Clearly, Assumption \ref{H2} holds for $\hat \a =\k -1$ and $\hat \b = 2$.
By Theorem \ref{Theorem4.1} and Proposition \ref{Cor3.10}, we have the following corollary.

\begin{corollary}[Convergence rate of TMil for A\"{i}t-Sahalia model]\label{Cor34}
Under the same assumption as Corollary \ref{Cor33},
the TMil scheme \eqref{PMil} for model \eqref{Ait1}
by setting
$$ R(\D)= L_1 \D^{- \frac{1}{ (2 \k-2 )\ve 4}}$$ has the property that
\begin{align}
\max_{0 \le k \le \lceil T/\D \rceil}\E \left [|X(\tk) - \pd (\hat X_\D (\tk))|^{2}\right ] & \le  C \D^{2}, \qu T >0. \label{Tem71611}
\end{align}
For any $p \ge 2$,
\begin{align*}
\max_{0 \le k \le \lceil T/\D \rceil}\E \left [| \pd (\hat X_\D (\tk))|^{p}\right ]
\ve
\max_{0 \le k \le \lceil T/\D \rceil}\E \left [| \pd (\hat X_\D (\tk))|^{-p}\right ]
< \infty.
\end{align*}

\end{corollary}
\begin{remark}\label{remark1}
We remark that the authors in \cite{Deng2023BIT} have derived the strong convergence rate of order arbitrarily close to $\frac{1}{2p}$ of the TEM method for AIT model \eqref{Ait1} for any $p \ge 1$ and $\eps \in (0, 1/2)$, i.e.,
\begin{align}\label{eq_724}
\max_{0 \le k \le \lceil T/\D \rceil}\E \left [|X(\tk) - \pd (X_\D (\tk))|^{p}\right ] & \le  C \D^{1/2-\eps},
\end{align}
see \cite[Theorem 3.1]{Deng2023BIT}. In contrast, we achieve the optimal strong convergence rate of order $\frac{1}{2}$.
\end{remark}
\subsection{CIR model}\label{cir_ap}
The Cox-Ingersoll-Ross (CIR) process is given by the SDE
\begin{align}\label{CIR_eq}
d X(t) = b_1( b_2 - X(t)) dt + \sigma \sqrt {X(t)} d B(t), \qu X(0)= X_0 >0,
\end{align}
where $b_1$, $b_2$, $\sigma$ are strictly positive constant parameters. Under the Feller condition
\begin{align}\label{Fell}
\varpi:= 2 b_1 b_2 / \sigma^2 >1,
\end{align}
 $X$ remains strictly positive almost surely. However Assumption \ref{Assu1} does not hold for CIR \eqref{CIR_eq}. Thus, our TEM can not apply to approximate CIR \eqref{CIR_eq} directly.
But, if we combine the Lamperti transform with our TEM, the strong $L^p$-convergence of order $1/2$ can be
derived for $p$ in a restricted parameter range.
\par
Applying the \Ito formula
to the Lamperti transform $Y= \sqrt{X}$ gives a  new SDE
\begin{align}\label{Cir_2}
dY(t) = f (Y(t)) dt - \frac{\sigma}{2}dB(t), \qu Y(0) = \sqrt{X_0},
\end{align}
where \begin{align*}
f(x) = \frac{\hat a}{x} + \hat b x, \qu \hat a = \frac{4 b_1 b_2 - \sigma^2}{8}, \qu \hat b := -\frac{b_1}{2}.
      \end{align*}
From \cite[p. 5]{Szpruch2012}, we have that
\begin{align*}
\sup_{0 \le t \le T} \E \left [ |X(t)|^ p\right ] < \infty, \qu \f p > - \varpi,
\end{align*}
and therefore,
\begin{align*}
\sup_{0 \le t \le T} \E \left [ |Y(t)|^ {-p}\right ] < \infty, \qu \f p  < \frac{4b_1 b_2}{\sigma^2} =2\varpi.
\end{align*}
Thus, for transformed SDE \eqref{Cir_2}, Assumptions \ref{H0}-\ref{H1} hold with $\a = 0$, $\b = 2$, $ p_1 < 2 \varpi$, $ p_0 =+ \infty$. Moreover, Assumption \ref{H2} also holds with $\hat \a = 1$, $\hat \b =2$.
According to  Proposition \ref{Cor3.1} and Theorem \ref{Theorem4.1}, we have the following results.

\begin{corollary}[Convergence rate of Lamperti TEM/TMil for CIR model]\label{Cor34}
Let $X$ be the CIR process \eqref{CIR_eq} and $\varpi >5 $, where $\varpi $ is  the parameter defined by \eqref{Fell}.
\begin{itemize}
  \item  Let
$Y_\D$ be
the TEM solution \eqref{PTEM} for \eqref{Cir_2}  by
setting
$ R(\D)= L_1 \D^{- \frac{1}{8}}$.
Then
for $p \in \left [1, \frac{\varpi}{5} \right ]$,
\begin{align}\label{Tem73111}
\max_{0 \le k \le \lceil T/\D \rceil}\E \left  [ \big | X (\tk) -  Y_\D (\tk)^2 \big | ^p \right ] & \le  C \D^{p/2}, \qu T >0.
\end{align}
  \item  Let
$\hat Y_\D$ be
the TMil solution \eqref{PMil} for \eqref{Cir_2}  by
setting
$ R(\D)= L_1 \D^{- \frac{1}{4}}$.
Then
\begin{align}\label{Tem731_2}
\max_{0 \le k \le \lceil T/\D \rceil}\E \left  [ \big | X (\tk) -  \hat Y_\D (\tk)^2 \big | \right ] & \le  C \D, \qu T >0.
\end{align}

\end{itemize}

\end{corollary}
\textbf{Proof.}
Let $Y$ be the solution of \eqref{Cir_2}.
In order to approximate the original CIR process, we observe that
\begin{align*}
X(\tk) -  Y_\D (\tk)^2=Y(\tk)^2 -  Y_\D (\tk)^2  = \Big ( Y (\tk) +  Y_\D (\tk) \Big )\Big ( Y (\tk) -  Y_\D (\tk) \Big ).
\end{align*}
Then the \Holder inequality gives
\begin{align*}
 & \E \left  [ \big | X (\tk) -  Y_\D (\tk)^2 \big | ^p \right ] = \E \left  [ \big | X (\tk) +  Y_\D (\tk) \big | ^p \big | X (\tk) -  Y_\D (\tk) \big | ^p \right  ]  \\
  & \le \left (  \E \left  [ \big | X (\tk) +  Y_\D (\tk) \big | ^{2p} \right ] \right )^{1/2}
 \left (  \E \left  [ \big | X (\tk) +  Y_\D (\tk) \big | ^{2p} \right ] \right )^{1/2}.
\end{align*}
  $ 1 \le p \le \frac{p_1}{10} < \frac{\varpi}{5}$ implies that  Condition \eqref{Cond7} is satisfied.
From \eqref{Moment_Y_1_1}, we obatian   $ \E \left  [  |  Y_\D (\tk)  | ^{2p} \right ]  < C$.
 Similarly,  $ \E \left  [  |  X_\D (\tk)  | ^{2p} \right ]  $ is finite. This, combined with \eqref{Eq626_9} lead to \eqref{Tem73111}.
 Similarly, by Theorem \ref{Theorem4.1} and  Corollary \ref{Cor3.10}, we can show that
 \eqref{Tem731_2}  also holds.
 $\Box$
 \begin{remark}\label{Rem4}
 Comparing Corollary \ref{Cor34} with Corollary 4.1 in \cite{FCJ2016SIAM}, where the authors
 proved the strong $L^p$-convegence of order $1/p$  of Lamperti projected EM for CIR, we observe that our Lamperti TEM has a significant improvement of the  $L^p$-convergence order   under slightly more stronger  condtion on the parameters.

 \end{remark}

\section{Numerical  experiments}\label{Sec5}
\label{sec:experiments}
In this section we compare  TEM \eqref{PTEM}  and TMil \eqref{PMil2} to some numerical schemes that we outline below
for the 3/2 and  the AIT models.
\begin{itemize}
  \item Euler-Maruyama scheme (EM) \cite{Deng2019PA}:
   $$Y_{k+1} = Y_k + f(Y_k) \D + g(|Y_k|)\D B_k, \qu Y_0 = X_0;$$
    \item Backward Euler-Maruyama scheme (BEM) \cite{Mao2011NA_Sahalia}:
    $$Y_{k+1} = Y_k + f(Y_{k+1}) \D + g(|Y_k|)\D B_k, \qu Y_0 = X_0;$$
\item Logarithmic Truncated EM scheme (log TEM) \cite{Tang_Mao2024ap}: 
$$Z_{k+1} = Z_k + F( \hat \pi (Z_{k+1})) \D + G(\hat \pi(Z_k))\D B_k, \qu Y_{k+1} = \textrm{e}^{Z_{k+1}}, \qu Z_0 = \log (X_0), \qu Y_0 = X_0,$$  where $F(x) = \textrm{e}^{-x} f(\textrm{e}^{x}) - 0.5 \textrm{e}^{-2x}|g(\textrm{e}^{x})|^2$, $G(x) = \textrm{e}^{-x} g(\textrm{e}^{x}) $, $\hat \pi (x) = (-R) \ve x \we R  $, $R = C + \frac{\log \D^{-1}}{\a \ve (\b +1)}$;
\item  Semi-implicit  Tamed EM  scheme   for AIT (STEM) \cite{WXJ2024na2}:
$$Y_{k+1} = Y_k + a_{-1}Y_{k+1}^{-1} \D + \left ( -a_0 + a_1 Y_k - \frac{a_2 Y_k^\k}{1+ \sqrt{\D} Y_k^\k} \right )\D +\frac{b Y_k^\theta}{1+ \sqrt{\D} Y_k^\k}\D B_k, \qu Y_0 = X_0, $$
\item Semi-implicit  Truncated  EM  scheme   for AIT (STEM2) \cite{WXJ2024na_Mil}:
$$Y_{k+1} = Y_k + a_{-1}Y_{k+1}^{-1} \D + \left ( -a_0 + a_1 Y_k - a_2 \check{Y}_k^\k \right )\D +
b \check{Y}_k^\theta\D B_k, \qu Y_0 = X_0, $$
where $\check{Y}_k =(-R) \ve x \we R$, $R= \D^{-1/(2\k -2)}$.
\end{itemize}
Denote by $Z^{\star}_j (T)$ the scheme $\star$ approximation at time $T$ and by 
$X^{Ref}_j(T)$ the reference solution calculated by the corresponding scheme $\star$ with small step size $\D=2^{-12}$, using the same Brownian motion path (the $j$th path).
The root mean square error (RMSE)  for the scheme $\star$ is defined by
\begin{align*}
\textrm{RMSE}^\star:=   \left ( \frac{1}{M} \sum_{j=1}^{M} |X^{Ref}_j(T) - Z^{\star}_j (T)|^2 \right )^{1/2},
\end{align*}
over  $M$ sample paths.
The strong error rates are computed by plotting RMSE against the step size on a log-log scale, and the strong rate of convergence is then retrieved using linear regression.
Moreover,   the event $ \mathbf {E}_\D$  is  defined  by
\begin{align}\label{eq810}
\mathbf {E}_\D: = \left \{\min _{ 0 \le k \le \lceil T/\D \rceil} X_\D (\tk) \le 0 \right \},
\end{align}
where $ X_\D$ is generated by the TEM scheme \eqref{PTEM2}.

\begin{figure}[!h]
  \centering
  \subfloat[RMSEs for 3/2 model with $\sigma =1$ ]{
    \label{fig1a} 
    \includegraphics[width=6.5cm,height=6.5cm]{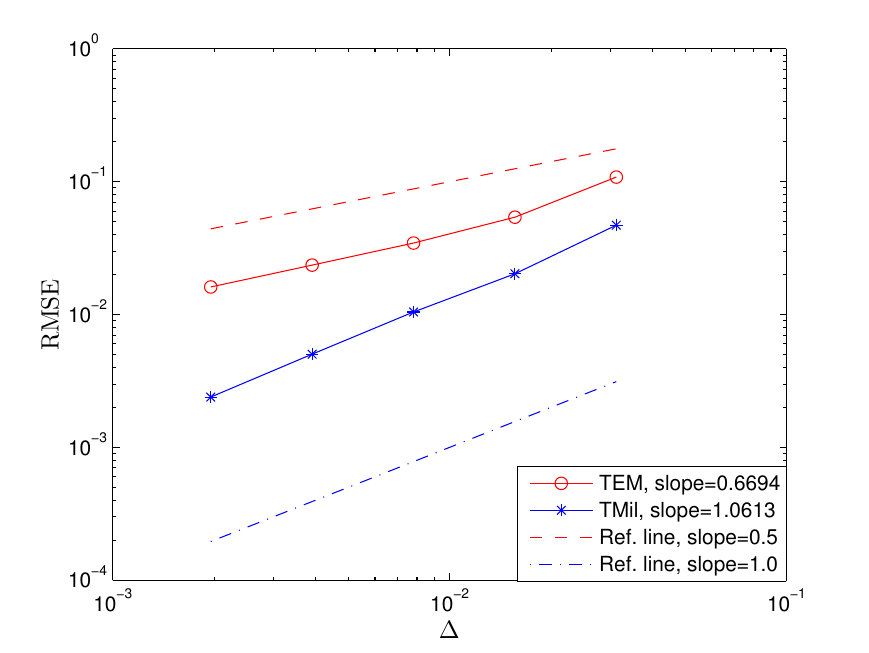}
  }
  \subfloat[RMSEs for AIT model]{
    \label{fig2a} 
    \includegraphics[width=6.5 cm,height=6.5cm]{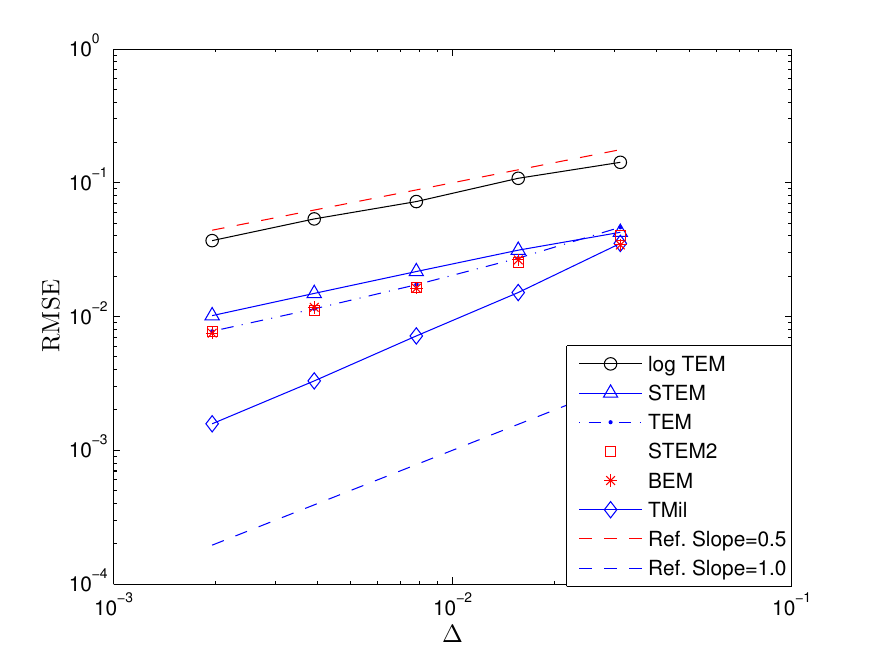}
  }
\caption{Convergence rates  for  3/2 and AIT models }
\end{figure}

\begin{example}\label{eg1}
Consider the following 3/2 model
\begin{align}\label{model61}
dX(t) = (4 X(t) - 4 X(t)^2 )dt+ \sigma X(t)^{3/2} dB(t), \qu X(0) =2,
\end{align}
where $\displaystyle \sigma ^2  <  \frac{4}{3}$.
 Now, we set $T=2$, $M=1000$, $ R(\D)= 50 \D^{- \frac{1}{2}}$. 
According to Corollaries \ref{cor31} and \ref{cor32}, TEM solution $Y_\D$ and TMil solution $\hat Y_\D $ has the property that 
\begin{align*}
\E |X(T) - Y_\D (T)|^2  \le C\D, \quad \textrm{and}\quad \E |X(T) - \hat Y_\D (T)|^2  \le C\D^2,
\end{align*}
respectively. RMSEs for different step sizes and rates
of TEM and TMil 
for 3/2 model are displayed in Table \ref{Tab3} and Figure \ref{fig1a}.
In the last line in Table \ref{Tab3}, we observe 
the empirical rates 
0.6138 and 1.0537 of TEM and TMil for $\sigma = 1/2$, 
0.6694 and 1.0613 of TEM and TMil for $\sigma = 1$,  
which shows that the observed  rates of convergence  are slightly higher than  the theoretical rates. 
\par
Moreover,   the fourth and the seventh column    in Table \ref{Tab3} 
lists respectively the
probabilities of
the TEM solutions  escaping from $\R_+$ for different step sizes under a small noise intensity with $\sigma =1/2$, and a higher intensity with $\sigma =1$.
We observe that as the step size goes to zero, this probability $\PP (\mathbf {E}_\D)$ tends to zero.
In other words,   for a sufficiently small step size,  say $\D = 2^{-7}$, the truncation $\pd $ from  \eqref{PTEM} does not execute for Example \ref{eg1}. In this context, our TEM  and the EM coincide with a large probability.
%
\end{example}
\begin{table}[htbp]
\caption{RMSEs and rates of TEM and TMil with different step sizes  for  3/2 model }\label{Tab3}
\begin{tabular*}{\textwidth}{@{\extracolsep\fill}lcccccc}
\toprule
& \multicolumn{3}{@{}c@{}}{$\sigma = 1/2$}
& \multicolumn{3}{@{}c@{}}{$\sigma = 1$} \\\cmidrule{2-4}\cmidrule{5-7}%
$\Delta$ & TEM & TMil & $\PP (\mathbf {E}_\D)$ & TEM & TMil &  $\PP (\mathbf {E}_\D)$ \\
\midrule
$2^{-3}$  &             &            & $0.03\%$  &          &             & $15.56\%$\\
$2^{-4}$  &             &            & $0.00\%$  &          &             & $2.90\%$\\
$2^{-5}$  & 2.2170e-02 & 1.3748e-02 & $0.00 \%$ & 1.0827e-01 & 4.7061e-02 & $0.20\%$\\
$2^{-6}$  & 1.3872e-02 & 6.3934e-03  & $0.00 \%$  & 5.4048e-02 & 2.0271e-02 & $0.01\%$\\
$2^{-7}$  & 8.6830e-03 & 3.1602e-03 & $0.00 \%$ & 3.4534e-02 & 1.0442e-02 & $0.00\%$\\
$2^{-8}$  & 5.9060e-03 & 1.4982e-03  & $0.00 \%$  &  2.3511e-02 & 5.0293e-03      & $0.00 \%$   \\
$2^{-9}$  & 4.0491e-03 & 7.3669e-04  & $0.00 \%$ & 1.6132e-02 &  2.3878e-03       & $0.00 \%$   \\
rate       &0.6138     & 1.0537 &  &0.6694      &1.0613             &     \\
\bottomrule  
\end{tabular*}
  \begin{tabular}{l}
  \end{tabular}
\footnotetext{Note: This is an example of table footnote. This is an example of table footnote this is an example of table footnote this is an example of~table footnote this is an example of table footnote.}
\end{table}
\begin{table}[htbp] 
\caption{Numerical results for different schemes for AIT }
\label{Tab5}
\begin{tabular*}{\textwidth}{@{\extracolsep\fill}ccccccc}
\toprule
 $\D$  & log TEM & STEM & STEM2 & BEM & TEM  & TMil\\
\midrule
$2^{-5}$ &  1.4201e-01 & 4.2682e-02  & 4.0054e-02 & 3.4322e-02 & 4.6424e-02 & 3.5164e-02\\
$2^{-6}$  & 1.0806e-01 & 3.1297e-02  & 2.5265e-02 & 2.6534e-02 &  2.7311e-02 & 1.5099e-02\\
$2^{-7}$ & 7.2431e-02 & 2.1658e-02   & 1.6477e-02 & 1.6354e-02 & 1.7300e-02 & 7.1460e-03\\
$2^{-8}$  & 5.3581e-02 & 1.4881e-02  & 1.1161e-02 & 1.1703e-02 & 1.1393e-02 & 3.2951e-03\\
$2^{-9}$ & 3.6888e-02 & 1.0145e-02 & 7.6605e-03 &  7.5196e-03 & 7.7554e-03 & 1.5746e-03\\
rate  & 0.4902 & 0.5218  &  0.5951 &    0.5566      &  0.6425& 1.1158\\
CPU time  & 4.86s & 4.19s  &  3.56s &  2161.61s        &  3.16s & 3.19s \\
\bottomrule  
\end{tabular*}
\footnotetext{Note: This is an example of table footnote. This is an example of table footnote this is an example of table footnote this is an example of~table footnote this is an example of table footnote.}
\end{table}

\begin{figure}[!h]
  \centering
    \includegraphics[width=14cm,height=7cm]{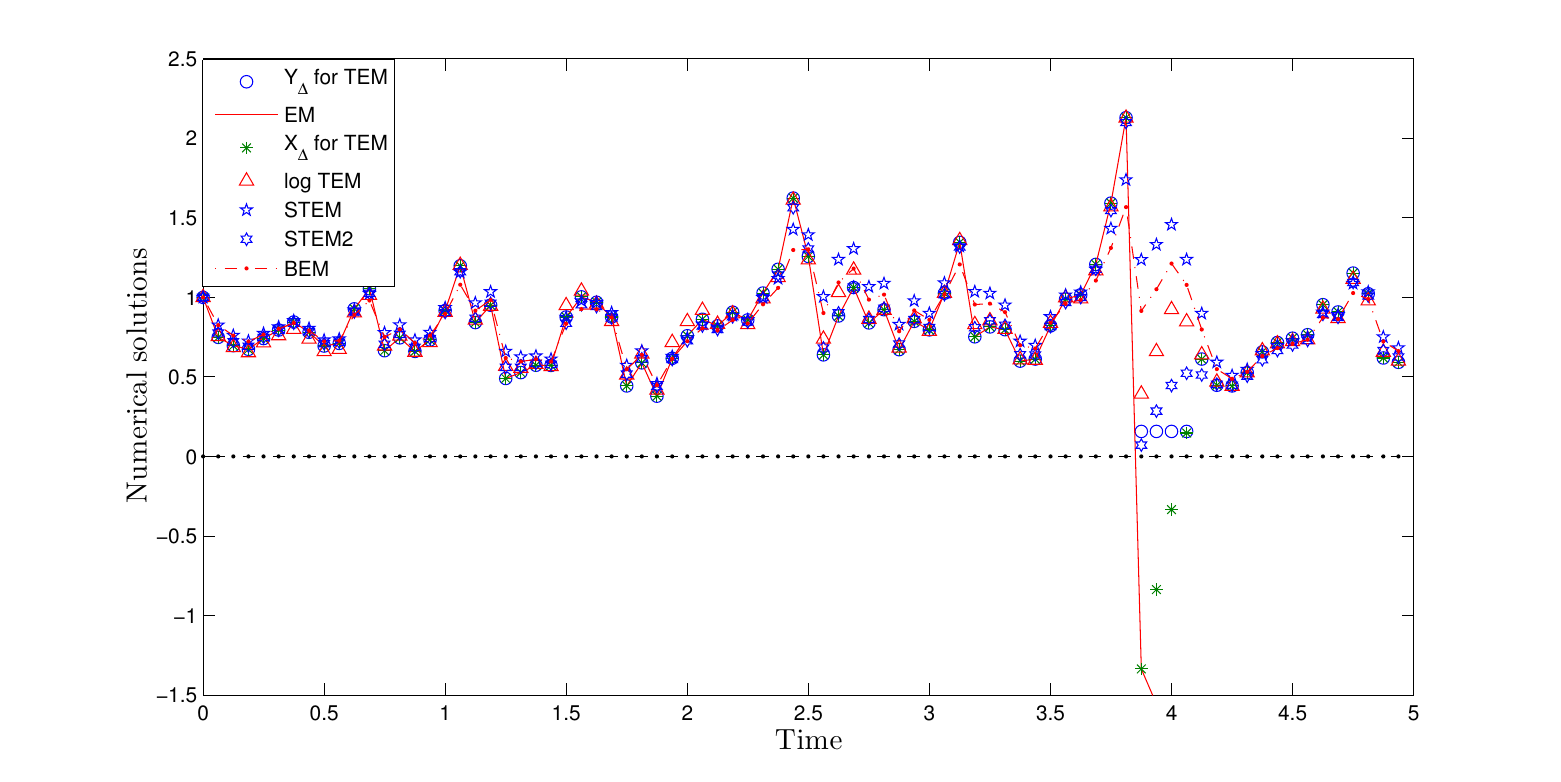}
\caption{Sample trajectories computed with  different schemes for  AIT }  \label{fig6}
\end{figure}

\begin{example}\label{ex2}
Consider the  AIT model \eqref{Ait1} with the following parameters
\begin{align*}
a_{-1}=1.5, \qu a_0 =2, \qu a_1=1, \qu a_2 = 2, \qu b =1, \qu \k =4, \qu \th =1.5, \qu X_0=1. 
\end{align*}
Clearly, the above setting satisfies the condition $\k +3 > 4 \th$. Thus, all the methods listed in Table \ref{Tab2} can be applied  to  this AIT. 
\par
Sample trajectories from  TEM, EM, log TEM, STEM, STEM2, and BEM  for AIT are plotted in Figure \ref{fig6}, which shows that  truncation occurs where 
the solution is close to zero and  thus 
these methods except EM
maintain positivity of the numerical approximations.

RMSEs and rates for different schemes are presented in Figure \ref{fig2a} and Table \ref{Tab5}.
In Figure \ref{fig2a}, we observe that 
log TEM has the largest error constant, whereas STEM, STEM2, BEM and TEM appear to have   smaller error constants.  
Furthermore, 
the graphs of the log TEM, the STEM, the STEM2, and the BEM 
seem parallel to the reference line with slope equal to $0.5$, while the TEM has a slightly bigger slope. 
Nevertheless, we observe in Figure \ref{fig2a} and Table \ref{Tab5} that TMil  reaches the optimal convergence rate $1$ and thus has a slight advantage over TEM when high accuracy is required. 
This behavior  is also confirmed  in the   
sixed line in Table \ref{Tab5}. As expected the empirical rates of convergence  are close to the theoretical ones.

Computational costs  as measured by CPU time in seconds   for difference methods are illustrated 
in the last line in Table \ref{Tab5}, where we set $M=10^4$, $T=2$ and  $\D = 2^{-12}$.
One clearly observes that the explicit  log TEM, STEM, STEM2, TEM and TMil  greatly
decreases the computational  time compared to the implicit BEM. 
Since in the implementation of the BEM, we have to solve numerically a non-linear equation at each time step.  This extra step brings questions about  the computing performance of BEM. 
%
%
%
%
%
%
%
%
%



\end{example}

\section{Appendix}\label{Sec6}
\subsection{Proof of Lemma \ref{Lem1}}\label{Append1}
\textbf{Proof.}
Write $R=R(\D)$ for simplicity.We first prove \eqref{Eq2}.
For $x,y \in (-\infty,R^{-1}]$ or $[R^{-1},R]$, or $[R,\infty)$, obviously \eqref{Eq2} holds.
Without loss of generality, we assume that $x \le y$.
For $R^{-1} \le x \le R \le y$, $\pd(x) = x$, $\pd(y) = R$, then
\begin{align}\label{AA1}
|\pd(x)- \pd(y)| = |x-R| =R -x \le y-x = |y-x|.
\end{align}
For $ x \le R^{-1} \le R \le y$, $\pd(x) = R^{-1} $, $\pd(y) = R$,
\begin{align}\label{AA2}
|\pd(x)- \pd(y)| = R - R^{-1} \le |y-x|.
\end{align}
For $ x \le R^{-1} \le y \le R$, $\pd(x) = R^{-1} $, $\pd(y) = y$,
\begin{align}\label{AA3}
|\pd(x)- \pd(y)| = y - R^{-1} \le |y-x|.
\end{align}
From \eqref{AA1}-\eqref{AA3}, we conclude that \eqref{Eq2} holds.

Now, we begin to prove \eqref{Eq11}.
For $x,y \in [R, \infty)$, or $ [R^{-1},R]$, or $(-\infty, R^{-1}]$, \eqref{Eq11} holds obviously. Without loss of generality, we assume that $x \le y$.
For $R^{-1} \le x \le R \le y $, $\pd(x) = x$, $\pd(y) = R$,
then by the \Holder inequality and Assumption \ref{Assu1}, we have
\begin{align}\label{Eq6}
& 2 \lan y-x, \fd (y) - \fd(x)\ran  + q_0 |\gd(y)-\gd(x)|^2  \\
& = 2\lan y-x, f(R)- f(x) \ran + q_0 |g(R)-g(x)|^2  \no \\
&  = 2 (y-x) \int^R_x f'(u)du + q_0  \Big | \int^R_x g'(u)du \Big |^2 \no \\
& \le  (y-R) \int^R_x 2 f'(u)du + (R-x) \int^R_x 2 f'(u)du + (R-x) \int^R_x q_0 |g'(u)|^2 du \no \\
& =  (y-R) \int^R_x 2 f'(u)du + (R-x)\int^R_x 2 \Big(  f'(u) + q_0 |g'(u)|^2 \Big ) du \no \\
& \le K (y-R) (R-x) + K |R-x|^2  \no \\
& \le 2K |y-x|^2. \no
\end{align}
Similarity, for $x \le R^{-1} < R \le y $, $\pd(x) = R^{-1}$, $\pd(y) = R$, we have
\begin{align}\label{Eq6_2}
& 2 \lan y-x, \fd (y) - \fd(x)\ran  + q_0 |\gd(y)-\gd(x)|^2 \no \\
& = 2\lan y-x, f(R)- f(R^{-1}) \ran + q_0 |g(R)-g(R^{-1})|^2  \no \\
&  = 2 (y-x) \int^R_{R^{-1}} f'(u)du + q_0  \Big | \int^R_{R^{-1}} g'(u)du \Big |^2 \no \\
& \le  (y-R + R^{-1} -x) \int^R_{R^{-1}} 2 f'(u)du + (R-R^{-1}) \int^R_{R^{-1}} 2 f'(u)du + (R-R^{-1}) \int^R_{R^{-1}} q_0 |g'(u)|^2 du \no \\
& =  (y-R + R^{-1} -x) \int^R_{R^{-1}} 2 f'(u)du + (R-R^{-1})\int^R_{R^{-1}} 2 \Big(  f'(u) + q_0 |g'(u)|^2 \Big ) du \no \\
& \le K (y-R + R^{-1} -x ) (R-x) + K |R-x|^2  \no \\
& \le 2K |y-x|^2.
\end{align}
For $x \le R^{-1} < y \le R $,  $\pd(x) = R^{-1}$, $\pd(y) = y$, we have
\begin{align*}
& 2 \lan y-x, \fd (y) - \fd(x)\ran  + q_0 |\gd(y)-\gd(x)|^2 \\
& = 2\lan y-x, f(y)- f(R^{-1}) \ran + q_0 |g(y)-g(R^{-1})|^2  \no \\
&  = 2 (y-x) \int^y_{R^{-1}} f'(u)du + q_0  \Big | \int^y_{R^{-1}} g'(u)du \Big |^2 \no \\
& \le  (R^{-1}-x) \int^y_{R^{-1}} 2 f'(u)du + (y-R^{-1})\int^y_{R^{-1}} 2 \Big(  f'(u) + q_0 |g'(u)|^2 \Big ) du \no \\
& \le K (R^{-1}-x) (y-R^{-1}) + K |y-R^{-1}|^2  \no \\
& \le 2K |y-x|^2. \no
\end{align*}
Thus, the proof is finished. $\Box$

\subsection{Proof of Lemma \ref{Lem3.3_2}} \label{Append2}

\textbf{Proof.}
For  $p \in \left [2,
\frac{p_0}{1+\a} \we  \frac{p_1}{\b}\right ]$ and $t \ge 0$, we conclude from Assumption \ref{H1} and \eqref{Eq39} that

\begin{align}\label{Eq151-9}
\E \left [ |f(X(t)) |^p \ve |g(X(t)) |^p \right]
& \le C \E \left [(1 + |X(t)|^{p (1+\a)} + |X(t)|^{-p\b} \right] \no  \\
& \le C (1 + \E [ |X(t)|^{p_0} ] + \E [|X(t)|^{-p_1}] ) \le C.
\end{align}
Moreover,
\begin{align*}
X(t+\D) - X(t)= \int_t^{t +\D} f(X(s))ds + \int_t^{t +\D} g(X(s))dB(s).
\end{align*}
Hence, by the \Ito isometry and \eqref{Eq151-9} as well as \eqref{Eq61}, we have
\begin{align}\label{Eq153}
& \E [| X(t+\D) - X(t)|^p] \\
& \le C_p \left ( \E \left [ \Big | \int_t^{t+\D} f(X(s))ds \Big |^p\right ]
+ \E \left [ \Big | \int_t^{t+\D} g(X(s))dB(s) \Big |^p \right ] \right ) \no \\
& \le  C_p \left ( \D^{p-1}\E \left [ \int_t^{t+\D}   \Big |  f(X(s)) \Big |^p ds  \right ]
+ \D^{p/2-1}\E \left [ \int_t^{t+\D}   \Big |  g(X(s)) \Big |^p ds  \right ]   \right ) \no \\
& \le C_p \left ( \D^{p-1} \int_t^{t+\D} \sup_{0 \le s \le T} \E \left [ \Big |  f(X(s)) \Big |^p  \right ] ds
+ \D^{p/2-1} \int_t^{t+\D} \sup_{0 \le s \le T} \E \left [ \Big |  g(X(s)) \Big |^p  \right ] ds   \right ) \no \\
& \le C_p \D^{p/2}. \no
\end{align}
For the  case of  $p \in (0,2)$, the required results follows by  the Lyapunov inequality.
$\Box$

\subsection{Proof of Lemma \ref{Lem3.1}} \label{Append3}
\textbf{Proof.}
Write $R = R(\D)$ for simplicity.
We conclude from
 \eqref{Eq42} that
\begin{align}\label{Eq82}
\pd (x) + \frac{1}{\pd(x)} \le |x|+ \frac{1}{|x|} + 1, \quad \f x \in \R_+.
\end{align}
Combining this with  \eqref{Eq32}, we observe that
\begin{align}\label{Eq84}
& | f(X(t)) - \fd(X(t))|^p \ve | g(X(t)) - \gd(X(t))|^p  \\
& \le C (1 + |X(t)|^{p\a} + |\pd (X(t))|^{p\a}
+ |X(t)|^{-p \b} + |\pd (X(t))|^{-p \b}
)  |X(t) - \pd (X(t))|^p \no \\
& \le C (1 + |X(t)|^{p\a} + |X(t)|^{-p \b})  |X(t) - \pd (X(t))|^p \no \\
& \le  C \frac{1}{R^p}  \Big ( 1+ \frac{1}{R^{p\a}} + \frac{1}{|X(t)|^{p \b}} \Big )\II_{ \{X(t) < R^{-1}\}}
+ C  |X(t)|^p  \Big ( 1+ |X(t)|^{p\a} + \frac{1}{R^{p \b}} \Big )\II_{\{ X(t) > R \}} \no \\
&
\le
\underbrace{  C \frac{1}{R^p}  \Big ( 1+ |X(t)|^{- p \b} \Big )\II_{ \{X(t) < R^{-1} \} }}_{=:\hat J_1}
+ \underbrace{C \Big ( 1+ |X(t)|^{p\a +p}\Big )\II_{ \{ X(t) > R\}}}_{=:\hat J_2}. \no
\end{align}
Consider first $\hat J_1$.
By the \Holder and Markov inequalities as well as Assumption \ref{H1}, we have
\begin{align}\label{Eq87}
\E \left [ |X(t)|^{- p \b} \II_{ \{X(t) < R^{-1} \} } \right ]
   & \le C \Big ( \E \left [ |X(t)|^{-p_1} \right ] \Big )^{p\b /p_1}
   \Big ( \PP \{ 1/X(t) > R \} \Big )^{1 - p \b /p_1} \no \\
   & \le C_{p_1} \Big ( \frac{ \E \left [ |X(t)|^{-p_1} \right ]}{R^{p_1}} \Big )^{1- p\b /p_1}
   \le  \frac{C_{p_1}}{R(\D)^{p_1 - p \b}}.
\end{align}
Thus,
\begin{align}\label{Eq91}
\E [ \hat J_1 ] & \le
 \frac{C}{R^p} \E   \Big [ |X(t)|^{- p \b} \II_{ \{X(t) < R^{-1} \} } \Big ]
 \le  \frac{C_{p_1}}{R(\D)^{p_1 - p \b +p }} = C_{p_1}\D^{\g ( p_1 - p \b +p]}.
\end{align}
Similarly, we also can show that
\begin{align}\label{Eq92}
\E [\hat J_2 ] & \le C  \E \left [   |X(t)|^{p\a +p}\II_{ \{ X(t) > R\}} \right ]
\no \\
& \le \frac{ C_{p_0}}{R(\D)^{p_0 - p( \a +1) }} = C_{p_0}\D^{\g [ p_0 - (p \a +p)]}.
\end{align}
By \eqref{Eq84}, \eqref{Eq91} and \eqref{Eq92}, we conclude that \eqref{Eq94} holds. \eqref{Eq94_3} can be obtained in the same way as the proofs of \eqref{Eq94}. The proof is finished.
$\Box$

\section*{Acknowledgment}
%
The authors would like to thank the
National Natural Science Foundation of China (12271003, 62273003,  72301173)
for their financial support.





\end{document}